\documentclass[reqno,oneside,12pt]{amsart}
\hoffset - 1.5cm
\usepackage{amsmath}
\usepackage{amssymb}
\usepackage{amstext}
\usepackage{amsopn}
\usepackage{amsfonts}

\newtheorem{Corollary}{Corollary}[section]
\newtheorem{Definition}{Definition}[section]
\newtheorem{Example}{Example}[section]
\newtheorem{Fact}{Fact}[section]
\newtheorem{Lemma}{Lemma}[section]
\newtheorem{Proposition}{Proposition}[section]
\newtheorem{Remark}{Remark}[section]
\newtheorem{Theorem}{Theorem}[section]
\newtheorem{Hipothesis}{Hipothesis}[section]

\newcommand{\ba}{\begin{array}}
\newcommand{\bc}{\begin{center}}
\newcommand{\bd}{\begin{description}}
\newcommand{\bdm}{\begin{displaymath}}
\newcommand{\be}{\begin{enumerate}}
\newcommand{\beq}{\begin{equation}}
\newcommand{\bdf}{\begin{Definition}}
\newcommand{\bex}{\begin{Example}}
\newcommand{\bft}{\begin{Fact}}
\newcommand{\bl}{\begin{Lemma}}
\newcommand{\bp}{\begin{Proposition}}
\newcommand{\br}{\begin{Remark}}
\newcommand{\bt}{\begin{Theorem}}
\newcommand{\bco}{\begin{Corollary}}
\newcommand{\bh}{\begin{Hipothesis}}
\newcommand{\ea}{\end{array}}
\newcommand{\ec}{\end{center}}
\newcommand{\ed}{\end{description}}
\newcommand{\edm}{\end{displaymath}}
\newcommand{\ee}{\end{enumerate}}
\newcommand{\eeq}{\end{equation}}
\newcommand{\edf}{\end{Definition}}
\newcommand{\eex}{\end{Example}}
\newcommand{\eft}{\end{Fact}}
\newcommand{\el}{\end{Lemma}}
\newcommand{\ep}{\end{Proposition}}
\newcommand{\er}{\end{Remark}}
\newcommand{\et}{\end{Theorem}}
\newcommand{\eco}{\end{Corollary}}
\newcommand{\eh}{\end{Hipothesis}}

\newcommand{\bH}{\mathbb{H}}
\newcommand{\bI}{\mathbb{I}}

\newcommand{\bK}{\mathbb{K}}

\newcommand{\bN}{\mathbb{N}}

\newcommand{\bR}{\mathbb{R}}

\newcommand{\bZ}{\mathbb{Z}}

\newcommand{\cC}{\mathcal{C}}

\newcommand{\cU}{\mathcal{U}}

\newcommand{\numsec}{\setcounter{Theorem}{0}\setcounter{Definition}{0}
\setcounter{Remark}{0} \setcounter{Lemma}{0} \setcounter{Fact}{0}
\setcounter{Proposition}{0} \setcounter{Corollary}{0}
\setcounter{Example}{0} \setcounter{equation}{0}
\setcounter{Property}{0}\renewcommand\theequation{\arabic{section}.\arabic{equation}}
\renewcommand\theTheorem{\arabic{section}.\arabic{Theorem}}
\renewcommand\theDefinition{\arabic{section}.\arabic{Definition}}
\renewcommand\theRemark{\arabic{section}.\arabic{Remark}}
\renewcommand\theLemma{\arabic{section}.\arabic{Lemma}}
\renewcommand\theFact{\arabic{section}.\arabic{Fact}}
\renewcommand\theProposition{\arabic{section}.\arabic{Proposition}}
\renewcommand\theCorollary{\arabic{section}.\arabic{Corollary}}
\renewcommand\theExample{\arabic{section}.\arabic{Example}}
\renewcommand\theProperty{\arabic{section}.\arabic{Property}}}

\newcommand{\numsubsec}{\setcounter{Theorem}{0}\setcounter{Definition}{0}
\setcounter{Remark}{0} \setcounter{Lemma}{0} \setcounter{Fact}{0}
\setcounter{Proposition}{0} \setcounter{Corollary}{0}
\setcounter{Example}{0} \setcounter{equation}{0}
\setcounter{Property}{0}\renewcommand\theequation{\arabic{section}.\arabic{subsection}.\arabic{equation}}
\renewcommand\theTheorem{\arabic{section}.\arabic{subsection}.\arabic{Theorem}}
\renewcommand\theDefinition{\arabic{section}.\arabic{subsection}.\arabic{Definition}}
\renewcommand\theRemark{\arabic{section}.\arabic{subsection}.\arabic{Remark}}
\renewcommand\theLemma{\arabic{section}.\arabic{subsection}.\arabic{Lemma}}
\renewcommand\theFact{\arabic{section}.\arabic{subsection}.\arabic{Fact}}
\renewcommand\theProposition{\arabic{section}.\arabic{subsection}.\arabic{Proposition}}
\renewcommand\theCorollary{\arabic{section}.\arabic{subsection}.\arabic{Corollary}}
\renewcommand\theExample{\arabic{section}.\arabic{subsection}.\arabic{Example}}
\renewcommand\theProperty{\arabic{section}.\arabic{subsection}.\arabic{Property}}}

\numberwithin{equation}{section} \errorcontextlines=0

\newcommand{\diag}{\mathrm{ diag \;}}

\newcommand{\sign}{\mathrm{ sign \;}}

\newcommand{\sone}{SO(2)}
\newcommand{\ds}{\displaystyle}
\newcommand{\nt}{\noindent}

\newcommand{\h}{\mathbb{H}}
\newcommand{\dg}{\nabla_{\sone}\mathrm{-deg}}
\newcommand{\ind}{\mathrm{ind}}

\begin{document}

\title[Autonomous Newtonian Systems]{Existence and  Continuation of Periodic Solutions \\  of Autonomous Newtonian Systems}

\author{Justyna Fura$^{\dag}$}
\address{
Faculty of Mathematics and Computer Science\\
Nicolaus Copernicus University \\
PL-87-100 Toru\'{n} \\ ul. Chopina 12/18 \\
Poland} \email{Justyna.Fura@mat.uni.torun.pl}

\author{Anna Ratajczak$^{\ddag}$}
\address{
Faculty of Mathematics and Computer Science\\
Nicolaus Copernicus University \\
PL-87-100 Toru\'{n} \\ ul. Chopina 12/18 \\
Poland} \email{aniar@mat.uni.torun.pl}

\author{S{\l}awomir Rybicki$^{\ddag}$}
\address{Faculty of Mathematics and Computer Science\\
Nicolaus Copernicus University \\
PL-87-100 Toru\'{n} \\ ul. Chopina 12/18 \\
Poland} \email{Slawomir.Rybicki@mat.uni.torun.pl}

\date{\today}

\thanks{$^{\dag}$ Research sponsored by the Doctoral Program in Mathematics at the Nicolaus Copernicus University,
Toru\'n, Poland}

\thanks{$^{\ddag}$ Partially supported by the Ministry of Scientific Research and Information Technology, Poland; under
grant  number 1 PO3A 009 27}

\begin{abstract}
In this article we study the existence and the continuation of periodic solutions of   autonomous Newtonian
systems. To prove the results we apply the infinite-dimensional version of the degree for $\sone$-equivariant
gradient operators defined by the third author in \cite{[RYB1]} and developed in \cite{[RYB2]}. Using the results
due to Rabier \cite{[RAB]} and Wang \cite{[WAN]} we show that the Leray-Schauder degree is not applicable in the
proofs of  our theorems, because it vanishes.
\end{abstract}
\maketitle

\numsec

\section{Introduction}

The first aim of this article is to study the existence of periodic solutions of the following system
\begin{equation}\label{introsys}
 \ddot x = - V'(x)
\end{equation}
where $V \in C^2(\bR^n,\bR)$ and $V'$ denotes the gradient of $V.$ We assume that $(V')^{-1}(0) =
\{p_1,\ldots,p_q\}$ is a finite set and that $V'(x)=V''(\infty) \cdot x + o(\|x\|)$ as $\|x\| \rightarrow \infty,$
where $V''(\infty)$ is a real symmetric $(n \times n)$-matrix.

Such a problem has been considered for $q=1$ by Amann and Zehnder, see \cite{[AMZE]}, and by Benci and Fortunato,
see \cite{[BEFO]}, for any $q \in \bN.$

Benci and Fortunato estimated the number of non-stationary $T$-periodic solutions of \eqref{introsys} as $T
\rightarrow \infty.$ To avoid some technicalities and to make the  proofs more transparent they assumed that all
the non-stationary $T$-periodic solutions are not $T$-resonant and that  potential $V$ is a Morse function. These
assumptions seems to be  restrictive and rather difficult to verify.

We relax these assumptions and therefore we obtain only the existence of at least one non-stationary $T$-periodic
solution of \eqref{introsys}. We formulate the sufficient conditions for the existence of non-stationary
$T$-periodic solutions of \eqref{introsys} in terms of $V''(p)$ and $\ind(-V',p),$ where $p \in
\{p_1,\ldots,p_q,\infty\}.$ It is worth to point out that we can treat problems with resonance at stationary
solutions and at the infinity. As a basic tool we use the degree   for $\sone$-equivariant gradient maps, see
\cite{[RYB1]}, \cite{[RYB2]}.

The second aim of this article is to study the continuation of non-stationary $T$-periodic solutions of the
following system
\begin{equation} \label{introfamisys}
  \ddot x = -V'_{\lambda}(x)
\end{equation}
where $V_{\lambda} \in C^2(\bR^n,\bR) \text{ for } \lambda \in \bR$  and potential $V_0$ possesses all the
properties of potential $V$ in  \eqref{introsys}. We formulate sufficient conditions for the existence of
connected sets of $T$-periodic solutions of \eqref{introfamisys} emanating from level $\lambda=0.$

We consider solutions of \eqref{introsys} and \eqref{introfamisys} as critical points of $\sone$-invariant
functional defined on a suitably chosen Hilbert space, which is an orthogonal representation of the group $\sone,$
see also \cite{[MAW-WIL],[RABIN]}. Gradient of this functional is an $\sone$-equivariant map in the form of a
compact perturbation of the identity.

It is known that the Conley index and the Morse theory are not suitable tools for the study of global bifurcations
and the continuation of critical points of functionals, see \cite{[AMB],[BOHME],[IZE0],[MARINO],[TAKENS]} for
discussion and examples. Since considered gradient is $\sone$-equivariant, the Leray-Schauder degree is not
applicable in our approach because it vanishes, see \cite{[RAB],[WAN]} and Remark \ref{rwc}. Therefore to prove
our results we apply the degree for $\sone$-e-\-qui\-va\-riant gradient maps. Degrees for $G$-equivariant gradient
maps has been defined in \cite{[Dancer0]} for $G=\sone.$ Next it was improved in \cite{[RYB1]}  and in
\cite{[GEB0]} for symmetries of any compact Lie group $G.$

After  introduction this article is organized in the following way.

In Section \ref{prelim}   we summarize without proofs the relevant material on the degree for $\sone$-equivariant
gradient maps. We finish this section with a continuation theorem of critical  orbits of $\sone$-invariant
functionals and the Rabier and Wang result concerning computation of the Leray-Schauder degree for
$\sone$-equivariant o\-pe\-ra\-tors.

The main results of  Section \ref{admhom} are Lemmas \ref{losp}, \ref{losinfty}. We construct in these lemmas
admissible $G$-equivariant gradient homotopies   for the class of  operators in the form of a compact perturbation
of a $G$-equivariant self-adjoint Fredholm operator $L.$

We use these homotopies in Section \ref{asl} in the case $L=Id \text{ and } G=\sone.$ Namely,  we simplify
computation of the degree for $\sone$-e\-qui\-va\-riant gradient maps at an isolated degenerate critical point
i.e. at a critical point with the isotropy group $\sone$ or at the infinity.

In Section \ref{persol} we formulate and prove the main results of this article. In Subsection \ref{persollin} we
study properties of the functional $\Phi_A$ associated to the linear system \eqref{LNS}. In Subsection
\ref{persolnonlin} we formulate and prove the sufficient conditions for the existence of non-stationary
$T$-periodic solutions of the nonlinear system \eqref{NNSO}, see Theorems \ref{istnienie}, \ref{istnienie2}. We
would like to emphasize that we cannot use in our approach the Leray-Schauder degree, see Remark \ref{rwc}. In
Subsection \ref{persolcon} we study continuation of non-stationary $T$-periodic solutions of family of nonlinear
equations \eqref{rodzina}, see Theorem \ref{twkontynuacja}.

In Section \ref{illu} we illustrate results proved in Section \ref{persol}.

\numsec
\section{Preliminaria}
\label{prelim}

In this section, for the convenience of the reader, we remind the main properties of the degree for
$\sone$-equivariant gradient maps defined in  \cite{[RYB1]}. This degree will be denoted briefly by
$\nabla_{\sone}\mathrm{-deg}.$ We finish this section with a theorem due to Rabier \cite{[RAB]} and Wang
\cite{[WAN]}.

Put $\displaystyle U(\sone)= \bZ \oplus \left( \bigoplus_{k=1}^{\infty} \bZ \right)$ and define actions
$$+ , \star :
U(\sone) \times U(\sone) \rightarrow U(\sone)$$
$$\cdot :
\bZ \times U(\sone) \rightarrow U(\sone)$$ as follows

\begin{align}
  \label{doda} \alpha + \beta =&\left(\alpha_0 + \beta_0, \alpha_1 + \beta_1,
 \ldots,\alpha_k+\beta_k,\ldots\right)   \\
 \label{m} \alpha \star \beta =&(\alpha_0 \cdot \beta_0, \alpha_0 \cdot \beta_1 +
\beta_0 \cdot \alpha_1,   \ldots, \alpha_0 \cdot \beta_k + \beta_0 \cdot \alpha_k, \ldots) \\
\gamma \cdot \alpha = &(\gamma \cdot \alpha_0, \gamma \cdot \alpha_1,   \ldots, \gamma \cdot \alpha_k, \ldots)
\end{align}
where $\alpha = (\alpha_0, \alpha_1, \ldots, \alpha_k, \ldots), \beta = (\beta_0, \beta_1, \ldots,\beta_k, \ldots)
\in U(\sone)$ and $\gamma \in \bZ.$ It is easy to check that $(U(\sone),+,\star)$ is a commutative ring with unit
$\bI=(1,0,\ldots) \in (\sone).$ Ring $(U(\sone),+,\star)$ is known as the tom Dieck ring of the group $\sone$. For
a definition of the tom Dieck ring $U(G),$ where $G$ is any compact Lie group, we refer the reader to
\cite{[DIECK]}.

\nt If $\delta_1,\ldots,\delta_q \in U(\sone),$ then we write $\ds \prod_{j=1}^q \delta_j$ for $\delta_1 \star
\ldots \star \delta_q.$ Moreover, it is understood that $\ds \prod_{j \in \emptyset} \delta_j= \bI \in U(\sone).$

\noindent Let $V$ be a real, finite-dimensional and orthogonal representation of the group $\sone.$ If $v \in V$
then  the subgroup $SO(2)_v=\{g \in \sone : g \cdot v =v\}$ is said to be the isotropy group of $v \in V.$ Let
$\Omega \subset V$ be an open, bounded and $\sone$-invariant subset and let $H \subset \sone$ be  closed subgroup.
Then we define
\begin{itemize}
  \item $\Omega^H=\{v \in \Omega : H \subset \sone_v\}=\{v \in \Omega : g v = v \: \forall \: g \in H\},$
  \item $\Omega_H=\{v \in \Omega : H = \sone_v\}.$
\end{itemize}

\nt Fix $k \in \bN$ and set $C^k_{\sone}(V,\bR) = \{f \in C^k(V,\bR) : f \text{ is } \sone-\text{invariant}\}.$
Let $f \in C^1_{\sone}(V,\bR).$ Since $V$ is an orthogonal representation, gradient  $\nabla f : V \rightarrow V$
is an $\sone$-e\-qui\-va\-riant $C^0$-map. If $H \subset \sone$ is a closed subgroup then $V^H$ is a
finite-dimensional representation of the group $\sone$ and  $\big(\nabla f\big)^H =\nabla \big(f_{\mid V^H }\big)
: V^H \rightarrow V^H$ is well-defined $\sone$-equivariant gradient map. Choose an open, bounded and
$\sone$-invariant subset $\Omega \subset V$ such that $(\nabla f)^{-1}(0) \cap
\partial \Omega = \emptyset.$ Under these assumptions we have defined in \cite{[RYB1]} the degree for
$\sone$-equivariant gradient maps $\nabla_{\sone}\mathrm{-deg}(\nabla f,\Omega) \in U(\sone)$ with coordinates
$$\nabla_{\sone}\mathrm{-deg}(\nabla f,\Omega)=$$$$=(\nabla_{\sone}\mathrm{-deg}_{\sone}(\nabla f,\Omega),
\nabla_{\sone}\mathrm{-deg}_{\bZ_1}(\nabla f,\Omega), \ldots, \nabla_{\sone}\mathrm{-deg}_{\bZ_k}(\nabla
f,\Omega), \ldots ).$$

\nt For $\gamma > 0$ and $v_0 \in V^{\sone}$ we put $D_{\gamma}(V,v_0) = \{v \in V :\ \mid v - v_0 \mid <
\gamma\}.$ In the following theorem we formulate the main properties of the degree for $\sone$-equivariant
gradient maps.

\bt[\cite{[RYB1]}]\label{wlas} Under the above assumptions the degree for $\sone$-e\-qui\-va\-riant gradient maps
has the following properties\et

\begin{enumerate}
  \item  \label{w1}   if
$\nabla_{\sone}\mathrm{-deg}(\nabla f, \Omega) \neq \Theta,$ then $(\nabla f)^{-1}(0) \cap \Omega \neq \emptyset,$
\item  if  $\nabla_{\sone}\mathrm{-deg}_H(\nabla f, \Omega) \neq 0,$ then $(\nabla f)^{-1}(0) \cap \Omega^H
\neq \emptyset,$
  \item  \label{w3} if
$\Omega = \Omega_0 \cup \Omega_1$ and $\Omega_0 \cap \Omega_1 = \emptyset,$ then
$$\nabla_{\sone}\mathrm{-deg}(\nabla f, \Omega) = \nabla_{\sone}\mathrm{-deg}(\nabla f, \Omega_0) +
\nabla_{\sone}\mathrm{-deg}(\nabla f, \Omega_1),$$
  \item  \label{w4}  if
$ \Omega_0 \subset \Omega$ is an open $\sone$-invariant subset and $(\nabla f)^{-1}(0) \cap \Omega \subset
\Omega_0,$ then
$$ DEG(\nabla f, \Omega) = DEG(\nabla f, \Omega_0),$$
\item   \label{w2} if    $ f \in C^1_{\sone}(V \times [0,1],\bR)$ is such that
$(\nabla_v f)^{-1}(0) \cap \left(\partial \Omega \times [0,1] \right) = \emptyset,$ then
$$\nabla_{\sone}\mathrm{-deg}(\nabla f_0, \Omega)=\nabla_{\sone}\mathrm{-deg}(\nabla f_1,\Omega),$$

\item \label{wzaw}  if $W$ is an orthogonal representation of the group  $\sone,$   then
  $$\nabla_{\sone}\mathrm{-deg}((\nabla f, Id), \Omega \times D_{\gamma}(W,0)) =
\nabla_{\sone}\mathrm{-deg}(\nabla f, \Omega),$$
\item  if  $ f \in C^2_{\sone}(V,\bR)$ is such that $\nabla f(0) = 0$ and $\nabla^2f(0)$
is an   $\sone$-equivariant self-adjoint isomorphism  then there is $\gamma > 0$ such that
$$\nabla_{\sone}\mathrm{-deg}(\nabla f, D_{\gamma}(V,0)) = \nabla_{\sone}\mathrm{-deg}(\nabla^2 f(0), D_{\gamma}(V,0)).$$
\end{enumerate}

\br \label{goodrem} Directly from the definition of the degree for $\sone$-equivariant gradient maps, see
\cite{[RYB1]}, it follows that
\begin{enumerate}
  \item  if $H \subset \sone$ is a closed subgroup and for any $v \in \Omega$  $\sone_v \neq H,$ then
  $\dg_H(\nabla f,\Omega)=0.$
  \item $$\dg_{\sone}(\nabla f,\Omega)=
  \begin{cases}
    0 & \text{ if } \:\: \Omega^{\sone}=\emptyset, \\
     1 & \text{ if } \:\: \Omega^{\sone}=\{0\} \subset V, \\
    \mathrm{deg}_B(\nabla f^{\sone}, \Omega^{\sone},0)  & \text{otherwise},
  \end{cases}
  $$
where $\mathrm{deg}_B$ denotes  the Brouwer degree.
\end{enumerate}
\er

\nt Below we formulate product formula for the degree for $\sone$-equivariant gradient maps.

\bt[\cite{[RYB2]}] \label{pft} Let $\Omega_i \subset V_i$ be an open, bounded and $\sone$-invariant subset of a
finite-dimensional, orthogonal representation $V_i$ of the group  $\sone,$ for $i=1,2.$ Let $f_i \in
C^1_{\sone}(V_i,\bR) $ be such that $\big(\nabla f_i \big)^{-1}(0) \cap
\partial \Omega_i = \emptyset,$  for $i=1,2.$  Then
$$
\nabla_{\sone}\mathrm{-deg}((\nabla f_1, \nabla f_2), \Omega_1 \times \Omega_2) =
\nabla_{\sone}\mathrm{-deg}(\nabla f_1, \Omega_1) \star \nabla_{\sone}\mathrm{-deg}(\nabla f_2, \Omega_2).
$$
\et

\noindent For $k \in \bN$ define a map $\rho^k : \sone \rightarrow  GL(2,\mathbb{R})$ as follows
\[
\rho^k (e^{i\cdot\theta})= \left[\begin{array}{lr}
\cos (k\cdot\theta)&-\sin (k\cdot\theta)\\
\sin (k\cdot\theta)&\cos (k\cdot\theta)
\end{array}\right]
\qquad 0\le\theta  <  2\cdot\pi.
\]
For $j,k \in \bN$ we denote by $\mathbb{R}[j,k]$ the direct sum of $j$ copies of $(\mathbb{R}^2 ,\rho^k)$, we also
denote by $\mathbb{R}[j,0]$  the trivial $j$-dimensional representation of $\sone$. We say that two
representations $V$ and $W$ are equivalent  if there exists an equivariant, linear isomorphism $T : V \rightarrow
W$. The following classic result gives a complete classification (up to equivalence) of finite-di\-men\-sio\-nal
representations of the group $\sone$ (see \cite{[ADM]}).

\begin{Theorem}[\cite{[ADM]}]
\label{tk} If $V$ is a finite-dimensional representation of $\sone$  then there exist finite sequences
$\{j_i\},\, \{k_i\}$ satisfying:\\ $ (*)\qquad k_i\in \{0\}\cup \bN,\quad  j_i\in \bN,\quad  1\le i\le r,  \: k_1
< k_2 < \dots  < k_r $
\\  such that $V$ is equivalent to $\displaystyle\bigoplus^r_{i=1} \mathbb{R}[j_i ,k_i]$. Moreover, the
equivalence class of $V$, ($V\approx\displaystyle\bigoplus^r_{i=1} \mathbb{R}[j_i,k_i]$) is uniquely determined by
$\{k_i\},\, \{j_i\}$ satisfying $(*)$.
\end{Theorem}

\nt We will denote by $m^-(L)$ the Morse index of a symmetric matrix $L.$

\nt To apply successfully any degree theory we need computational formulas for this invariant. Below we show how
to compute degree for $\sone$-equivariant gradient maps of a linear, self-adjoint, $\sone$-equivariant
isomorphism.

\bl[\cite{[RYB1]}] \label{lindeg} If $V \approx \bR[j_0,0] \oplus \bR[j_1,k_1] \oplus \ldots \oplus \bR[j_r,k_r],$
 $L : V \rightarrow V$ is a self-adjoint, $\sone$-equivariant, linear isomorphism and $\gamma > 0$ then \el
\begin{enumerate}
  \item $L= \diag (L_0,L_1,\ldots,L_r),$
  \item $$\nabla_{\sone}\mathrm{-deg}_H(L,D_{\gamma}(V,0))=
 \begin{cases}
    (-1)^{m^-(L_0)}, & \text{ for } H = \sone, \\
   \displaystyle (-1)^{m^-(L_0)} \cdot  \frac{m^-(L_i)}{2}, & \text{ for } H = \bZ_{k_i}\\
     0, & \text{ for } H \notin \{\sone, \bZ_{k_1}, \ldots, \bZ_{k_r}\},
  \end{cases}
$$
  \item in particular, if $L=-Id,$ then
 $$\nabla_{\sone}\mathrm{-deg}_H(-Id,D_{\gamma}(V,0))=
 \begin{cases}
    (-1)^{j_0}, & \text{ for } H = \sone, \\
   \displaystyle (-1)^{j_0}  \cdot j_i, & \text{ for } H = \bZ_{k_i},\\
     0, & \text{ for } H \notin \{\sone, \bZ_{k_1}, \ldots, \bZ_{k_r}\}.
  \end{cases}
$$
\end{enumerate}

\noindent Let $(\bH, \langle \cdot,\cdot \rangle _{\bH})$ be an infinite-dimensional, separable Hilbert space
which is an or\-tho\-go\-nal representation of the group $\sone$ and let $C_{\sone}^1(\bH ,\bR)$ denote the set of
$\sone$-invariant $C^1$-functionals. Fix $\Phi \in C_{\sone}^1(\bH ,\bR)$   such that
\begin{equation}
\label{odwzfun} \nabla \Phi(u)=   u -  \nabla \eta(u),
\end{equation}
where  $\nabla  \eta : \bH   \rightarrow \bH$ is an $\sone$-equivariant   compact operator. Let $\cU \subset \h$
be an open, bounded and $\sone$-invariant set such that $\left(\nabla \Phi \right)^{-1}(0) \cap \partial \cU =
\emptyset.$ In this situation $\ds \dg(Id - \nabla \eta, \cU) \in U(\sone)$ is well-defined, see \cite{[RYB1]} for
details  and properties of this degree.

\nt The following theorem is an infinite-dimensional generalization of Theorem \ref{pft}.

\bt \label{pinf} Let $\cU_i \subset \h_i$, be an open, bounded and $\sone$-invariant subset of
infinite-dimensional, orthogonal representation $\h_i$ of the group $\sone,$ for $i=1,2.$ Let $f_i \in
C^1_{\sone}(\h_i,\bR)$ be such that
\begin{enumerate}
  \item $\nabla f_i = Id -\nabla \eta_i$ is an operator in the form of a compact perturbation of the identity, for $i=1,2,$
  \item $\big(\nabla f_i \big)^{-1}(0) \cap \partial \cU_i =\emptyset,$ for $i=1,2.$
\end{enumerate}
Then the following formula holds true

$$\dg((Id - \nabla \eta_1,Id - \nabla \eta_2),\cU_1 \times \cU_2)=$$ $$=
\dg(Id - \nabla \eta_1,\cU_1) \star \dg(Id - \nabla \eta_2,   \cU_2).$$
 \et

\nt Let $L : \h \rightarrow \h$ be a linear, bounded, self-adjoint, $\sone$-equivariant operator with spec\-trum
$\sigma(L)=\{\lambda_i\}.$ By $V_L(\lambda_i)$ we will denote eigenspace of $L$ corresponding to the
e\-igen\-va\-lue $\lambda_i$ and we put $\mu_L (\lambda_i)=\dim V_L(\lambda_i).$ In other words $\mu_L(\lambda_i)$
is the multiplicity of the eigenvalue $\lambda_i.$ Since operator $L$ is linear, bounded, self-adjoint, and
$\sone$-equivariant, $V_L(\lambda_i)$ is a finite-dimensional, orthogonal representation of the group $\sone.$ For
$\gamma > 0$ and $v_0 \in \h^{\sone}$ put $D_{\gamma}(\h,v_0) = \{v \in \h :\  \mid v - v_0 \mid < \gamma\}.$

\nt Combining Theorem 4.5 in  \cite{[RYB1]} with Theorem \ref{pft} we obtain the following theorem.

\bt  \label{dizom}  Under the above assumptions if $1 \notin \sigma(L),$ then $$\dg(Id - L, D_{\gamma}(\h,0))=
\prod_{\lambda_i > 1} \dg(-Id, D_{\gamma}(V_L(\lambda_i),0)) \in U(\sone).$$ It is understood that if $\sigma(L)
\cap [1,+\infty)=\emptyset,$ then
$$\dg(Id - L, D_{\gamma}(\h,0))=\bI \in U(\sone).$$
\et

\nt Below we formulate the continuation theorem for   $\sone$-equivariant gradient o\-pe\-ra\-tors in the form of
a  compact perturbation of the identity. In other words we study continuation of critical orbits of
$\sone$-invariant $C^1$-functionals. The proof   this theorem is standard, but in this proof we have to replace
the Leray-Schauder degree with the degree for  $\sone$-equivariant gradient operators.

\bt \label{abscont} Let $\Phi \in C_{\sone}^1(\bH \times \bR,\bR)$  be such that $\ds \nabla_u \Phi(u,\lambda)= u
- \nabla_u \eta(u,\lambda),$ where $\nabla  \eta : \bH  \times \bR \rightarrow \bH$ is an $\sone$-equivariant
compact operator. Fix an open, bounded and $\sone$-invariant subset $\cU \subset \h$ and $\lambda_0 \in \bR$ such
that

\begin{enumerate}
  \item $\left(\nabla_u \Phi(\cdot,\lambda_0 )\right)^{-1}(0) \cap \partial \cU = \emptyset,$
  \item $\ds \nabla_{\sone}-\mathrm{deg}(\nabla_u \Phi(\cdot,\lambda_0),\cU) \neq \Theta \in U(\sone).$
\end{enumerate}
Then there exists continua (closed connected sets) $\cC^{\pm} \subset \h \times \bR,$ with
$$\cC^- \subset \left((-\infty,\lambda_0] \times \h\right) \cap \left(\nabla_u \Phi(\cdot,\lambda_0 )\right)^{-1}(0),$$
$$\cC^+ \subset \left([\lambda_0,+\infty) \times \h\right) \cap \left(\nabla_u \Phi(\cdot,\lambda_0 )\right)^{-1}(0),$$
and for both $\cC=\cC^{\pm}$ the following statements are valid
\begin{enumerate}
  \item $\cC \cap (\{\lambda_0\} \times \cU) \neq \emptyset,$
  \item either $\cC$ is unbounded or else $\cC \cap (\h \setminus cl(\cU )) \neq \emptyset.$
\end{enumerate}

\et

\nt We finish this section with a special case, $G=\sone,$ of the theorem due to  Rabier and Wang, see
\cite{[RAB],[WAN]}.

\bt \label{rawa} Let $\cU \subset \h$ be an open, bounded and $\sone$-invariant set and let $f \in
C^0_{\sone}(cl(\cU),\h)$ be an operator in the form   of a compact perturbation of the identity such that $0
\notin f(\partial \cU).$ Then $\ds \mathrm{deg}_{\rm LS}(f,\cU,0)=\mathrm{deg}_{\rm LS}(f^{\sone},\cU^{\sone},0),$
where $\mathrm{deg}_{\rm LS}$ denotes the Leray-Schauder degree.\et

\numsec

\section{Admissible G-equivariant Gradient Homotopies}
\label{admhom}

This section is of technical nature. We prove here the splitting lemma at the origin and at the infinity. In fact
we construct admissible $G$-equivariant gradient homotopies. In the next sections, using these homotopies, we will
compute the degree for $\sone$-equivariant gradient maps.

\nt Let $G$ be any compact Lie group and let $\Phi \in C^{2}_{G}(\h,\mathbb{R})$ has the following form
\begin{equation} \label{funkcjonal}
 \Phi(x) =\displaystyle{
\frac{1}{2}\langle Lx,x \rangle_{\h} - g(x) },
\end{equation}
where $\nabla g:\h \rightarrow \h$ is a $G$-equivariant compact operator and
\begin{enumerate}
\item[(F.1)] $L:\h \rightarrow \h$ is a $G$-equivariant self-adjoint Fredholm  operator.
\end{enumerate}
Assume that for $p \in \{0, \infty\}$
\begin{equation}\label{postacfunkcjonalu}
 \Phi(x) =\displaystyle{
\frac{1}{2}\langle (L-L_p)x,x \rangle_{\h} + \eta_p(x) },
\end{equation}
i.e.
\begin{equation*}
 \Phi(x) =\displaystyle{
\frac{1}{2}\langle (L-L_0)x,x \rangle_{\h} + \eta_0(x) =\frac{1}{2}\langle (L-L_{\infty})x,x \rangle_{\h} +
\eta_{\infty}(x)},
\end{equation*}
 where
\begin{enumerate}
\item[(F.2)] $\nabla \eta_p : \h \rightarrow \h$ is a $G$-equivariant, compact operator,
\item[(F.3)] $L_p : \h \rightarrow \h$ is a linear, $G$-equivariant,
self-adjoint and compact operator,
\item[(F.4)] $\| \nabla^2 \eta_p(x) \|  \rightarrow 0$ as $\| x \| \rightarrow p,$
\item[(F.5)] $0 \in \sigma(L-L_p),$
\item[(F.6)] $p$ is an isolated critical point of $\Phi$.
\end{enumerate}

\nt We treat $p=\infty$ as a critical point of $\Phi$ with Hessian $\nabla^2\Phi(\infty)=L-L_{\infty}.$ Moreover,
we say that $\infty$ is an isolated critical point if $\big(\nabla \Phi\big)^{-1}(0)$ is bounded.

\nt We will denote by $V_p$ and $W_p$ the kernel and the image of $\nabla^{2}\Phi(p) = L-L_p$, respectively.
Notice that $V_p$ and $W_p$ are finite and infinite-dimensional orthogonal representation of the group G,
respectively. Since the operator $L-L_p$ is self-adjoint, $\h = V_p \oplus W_p$. Put $A_p=(L-L_p)_{|W_p} : W_p
\rightarrow W_p.$ Notice that the operator $A_p$ is an isomorphism. From now on $\pi_p: \h \rightarrow W_p$ and
$Id-\pi_p: \h \rightarrow V_p$ stand for $G$-equivariant, orthogonal projections. Set $\Phi_1^p = (Id-\pi_p) \circ
\nabla \Phi$ and $\Phi_2^p = \pi_p \circ \nabla \Phi$.

\nt The following two versions of the implicit function theorem   will allow us to construct admissible
$G$-equivariant gradient homotopies.

\bt  \label{uwiklp}  Let $\Phi \in C^1_{G}(\h, \mathbb{R})$ be a functional given by \eqref{funkcjonal}. Suppose
that  $\Phi(x) = \displaystyle{\frac{1}{2}\langle (L-L_0)x,x \rangle_{\h}} + \eta_0(x)$ and that assumptions
$(F.1)-(F.5)$ are satisfied for $p=0$. Then there exist $\varepsilon > 0$ and $G$-equivariant, $C^{1}$-mapping
$w_0 : D_{\varepsilon}(V_0,0) \rightarrow W_0$ such that
\begin{enumerate}
    \item[\rm(i)] $w_0(0)=0, D w_0(0)=0,$
    \item[\rm(ii)] $\Phi_2^0(v,w)=0$ for $v \in D_{\varepsilon}(V_0,0)$ iff $w=w_0(v).$
\end{enumerate}
\et
\begin{proof}
The existence of $w_0 : D_{\varepsilon}(V_0,0) \rightarrow W_0$ we obtain from the nonequivariant version of the
implicit function theorem, where $\varepsilon > 0$ is   sufficiently small. What is left is to show that $w_0$ is
$G$-equivariant. Since the operator $\Phi_2^0$ is $G$-equivariant, $\Phi_2^0(gv, gw_0(v)) =0$  for all $v \in
D_{\varepsilon}(V_0,0), g \in G.$ Moreover, we have and $\Phi_2^0(gv, w_0(gv))=0$ for all $v \in
D_{\varepsilon}(V_0,0), g \in G.$ From the uniqueness of $w_0$ we obtain $gw_0(v) = w_0(gv)$ for all $v \in
D_{\varepsilon}(V_0,0)$ and $g \in G,$ which completes the proof.
\end{proof}

\bt \label{uwiklinfty} Let $\Phi \in C^1_{G}(\h, \mathbb{R})$ be a functional given by \eqref{funkcjonal}. Suppose
that $\ds\Phi(x)=\frac12 \langle (L- L_{\infty})x, x \rangle_{\h} + \eta_{\infty}(x)$ and that assumptions
$(F.1)-(F.5)$ are sa\-ti\-sfied for $p= \infty$. Then there is $\beta_0 > 0$ and $G$-equivariant $C^1$-mapping
$w_{\infty}: V_{\infty} \setminus cl(D_{\beta_0}(V_{\infty},0)) \rightarrow W_{\infty}$ such that
$\Phi_2^{\infty}(v,w) = 0$ for $v \in V_{\infty} \setminus cl(D_{\beta_0}(V_{\infty},0))$ iff $w= w_{\infty}(v).$
\et

\nt The above theorem has been proved as a part of the proof of Lemma 4.3 in \cite{[Bartsch-Li]}.
$G$-e\-qui\-va\-riance follows in the same way  as in the proof of  Theorem \ref{uwiklp}.

\nt As a consequence of Theorems \ref{uwiklp}, \ref{uwiklinfty} we obtain the following corollary.

\bco \label{cpi} Let $\Phi \in C^1_{G}(\h, \mathbb{R})$ be a functional given by \eqref{funkcjonal}. Fix $p \in
\{0,\infty\}$ and assume that  functional  $\ds\Phi(x)=\frac12 \langle (L- L_p)x, x \rangle_{\h} + \eta_p(x)$
satisfies assumptions $(F.1)-(F.5).$ If $\{x_n\} \subset (\nabla \Phi)^{-1}(0)$   converges to $p$ then there is
$n_0 \in \bN$ such that for any $n > n_0$ there is $v_n \in V_p$ such that $G_{x_n}=G_{v_n}.$\eco
\begin{proof} Fix $p=0.$ From Theorem \ref{uwiklp} it follows that there is $\varepsilon > 0$ such that if
$x_n \in (\nabla \Phi)^{-1}(0) \cap D_{\varepsilon}(\bH,p)$ then $x_n=(v_n,w_p(v_n)) \in \bH,$ where $w_p :
D_{\varepsilon}(V_p,p) \rightarrow W_p$ is a $G$-e\-qui\-va\-riant map. Since $w_p$ is $G$-e\-qui\-va\-riant,
$G_{v_n} \subset G_{w_p(v_n)}$ for any $v_n \in D_{\varepsilon}(V_p,p).$ Consequently,
$G_{x_n}=G_{(v_n,w_p(v_n))}=G_{v_n} \cap G_{w_p(v_n)}=G_{v_n},$ which completes the proof. The same proof remains
valid for $p=\infty,$ but we have to replace Theorem \ref{uwiklp} with Theorem \ref{uwiklinfty}.
\end{proof}

\nt Let us consider continuous, $G$-equivariant  extensions $\tilde{w}_p : V_p \rightarrow W_p$ of $w_p, p=0,
\infty,$ defined in Theorems \ref{uwiklp}, \ref{uwiklinfty}. To shorten notation we continue to write $w_p$ for
$\tilde{w}_p, p=0, \infty.$

\nt For $p \in \{0, \infty\}$ consider a family of $G$-invariant $C^{1}$-functionals $H_p \in C_G^1(\h \times
[0,1], \mathbb{R})$, defined in a following way \beq \label{homotopia}
\begin{split}
   H_p((v,w),t)  &=  \frac{1}{2}\langle A_p(w),w\rangle_{\h} +
  \frac{1}{2}t(2-t) \langle A_p(w_p(v)),w_p(v)\rangle_{\h} +  \\
    &+ t \eta_p(v,w_p(v))+
  (1-t)\eta_p(v,w+t w_p(v)),
\end{split}
\eeq where $(v,w)\in \h=V_p\oplus W_p.$

\nt This family has been introduced by Dancer in  \cite{[Dancer]}. Recall that homotopy $H : \Omega \times [0,1]
\rightarrow \mathbb{R},$ where $\Omega$ is an open, bounded, $G$-invariant subset of an orthogonal
$G$-representation $\h$, is said to be $\Omega$-admissible iff
\begin{enumerate}
  \item[(2)] $\nabla H(\cdot,t) : (cl(\Omega),\partial \Omega) \rightarrow (\h,\h\backslash \{0\})$
  is a gradient, $G$-equivariant mapping for all $t \in [0,1],$
  \item[(3)] $\nabla H(\cdot,t) = L - \nabla g_t,$ where $\nabla g(\cdot,t)=\nabla g_t$ and $\nabla g: \h \times [0,1] \rightarrow
\h $ is a compact operator.
\end{enumerate}
By $\nabla H$ we denote the gradient of $H$ with respect to the coordinate $x \in \h$.

\bl\label{postac} For $p \in \{ 0, \infty\}$ the family $H_p \in C_G^1(\h \times [0,1], \mathbb{R})$ given
 by the formula \eqref{homotopia}
is of the form $\nabla H_p(\cdot,t) = L-\nabla g_t$ for $t\in[0,1],$ where $\nabla g: \h \times [0,1] \rightarrow
\h $ is a compact operator. \el

\begin{proof}
\noindent Observe that
\begin{align*}
\nabla H_p((v,w),t) &= A_p(w)+t(2-t)[D w_p(v)]^{T} A_p(w_p(v))
+t(Id-\pi_p)\nabla\eta_p(v,w_p(v))+\\& +t[D
w_p(v)]^{T}\pi_p\nabla\eta_p(v,w_p(v))
+(1-t)\nabla\eta_p(v,w+tw_p(v))+\\& + (1-t)t[D
w_p(v)]^{T}\pi_p\nabla\eta_p(v,w+t w_p(v)).
\end{align*}
\noindent From definition we have $A_p(w)=(L-L_p)(v,w).$ Recall that $L_p$ and $\nabla \eta_p$ are compact
operators and $V_p$ is finite dimensional space. To complete the proof we use following facts
\begin{enumerate}
\item[(i)] superposition of compact and continuous mappings is compact,
\item[(ii)] continuous, finite-dimensional mapping is compact,
\item[(iii)] continuous mapping defined on a finite dimensional Banach space is compact.
\end{enumerate}
\end{proof}

\nt We finish this section with splitting lemmas at the origin and at the infinity.


\bl (Splitting lemma at the origin) \label{losp} Suppose that functional $\Phi \in C^{2}_{G}(\h,\mathbb{R})$ is
given by  formula \eqref{funkcjonal}. Assume additionally that  for $p=0$  there is representation $\Phi(x) =
\frac{1}{2}\langle(L-L_0)x,x\rangle_{\h} + \eta_0(x) $ and assumptions (F.1)-(F.6) hold. Then there exist
$\alpha_0 > 0$ and $G$-equivariant gradient homotopy $\nabla H_0 : (V_0 \oplus W_0) \times [0,1] \rightarrow \h,$
satisfying the following conditions
\begin{enumerate}
  \item[\rm (1)] $\nabla H_0^{-1}(0) \cap (cl(D_{\alpha_0}(V_0,0)) \times cl(D_{\alpha_0}(W_0,0))
  \times [0,1]) = \{0\}\times [0,1],$
  \item[\rm (2)] $\nabla H_0((v,w),t)=(L-\nabla g_t)(v,w), $ for $t \in [0,1], (v,w)
\in V_0\oplus W_0,$ where $\nabla g_t=\nabla g(\cdot,t)$ and $\nabla g : \h \times [0,1] \rightarrow \h$ is a
compact mapping,
  \item[\rm (3)] $\nabla H_0((v,w),0) = \nabla \Phi (v,w)$,
  \item[\rm (4)] there exists a $G$-equivariant,  gradient mapping
  $\nabla\varphi_0 : (V_0,0) \rightarrow (V_0,0)$ such that
  $\nabla H_0((v,w),1) = (\nabla\varphi_0(v),A_0(w))= (\nabla \varphi_0(v),(L-L_0)_{|W_0}(w)).$
\end{enumerate}
\el
\begin{proof}
\noindent Let $H_0 : (V_0 \oplus W_0) \times [0,1] \rightarrow \mathbb{R}$ will be defined by formula
\eqref{homotopia} for $p=0$, i.e.
\begin{align*}
   H_0((v,w),t)  &= \frac{1}{2}\langle A_0(w),w\rangle_{\h} +
  \frac{1}{2}t(2-t) \langle A_0(w_0(v)),w_0(v)\rangle_{\h} +  \\
    &+ t \eta_0(v,w_0(v))+
  (1-t)\eta_0(v,w+t w_0(v)).
\end{align*}
\noindent Clearly $\nabla H_0 : (V_0\oplus W_0) \times [0,1] \rightarrow \h$ is a gradient, $G$-equivariant
homotopy. Condition $(1)$ has been verified for example in \cite{[Krysz-Szul]}, for $\alpha_0 > 0$ taken as a
sufficiently small number. Lemma \ref{postac} yields $(2)$. To complete the proof notice that
\begin{enumerate}
    \item[\rm (1)] $\nabla H_0((v,w),0) = \nabla \Phi(v,w) $,
    \item[\rm (2)] $\nabla H_0((v,w),1) = (\nabla\varphi_0(v),A_0(w)),$ where
    $\nabla \varphi_0(v) = \Phi_1^0 (v,w_0(v)).$
\end{enumerate}
\end{proof}

\bl (Splitting lemma at the infinity) \label{losinfty}  Suppose that functional $\Phi \in
C^{2}_{G}(\h,\mathbb{R})$ is  given by   formula \eqref{funkcjonal}. Assume additionally that for $p=\infty$ there
is representation $\Phi(x) = \frac{1}{2}\langle(L-L_{\infty})x,x\rangle_{\h} + \eta_{\infty}(x) $ and assumptions
(F.1)-(F.6) hold. Then there exist number $\alpha_{\infty} > 0$ and $G$-equivariant gradient homotopy $\nabla
H_{\infty}: (V_{\infty} \oplus W_{\infty}) \times [0,1] \rightarrow \h,$ satisfying the following conditions
 \be
\item[(1)] $\nabla H_{\infty}^{-1}(0) \subset cl(D_{\alpha_{\infty}}(V_{\infty},0)) \times
cl(D_{\alpha_{\infty}}(W_{\infty},0)) \times[0,1],$
\item[ (2)] $\nabla H_{\infty}((v,w),t)=(L-\nabla g_t)(v,w), $ for $t \in [0,1], (v,w)
\in V_{\infty} \oplus W_{\infty},$ where $\nabla g_t= \nabla g(\cdot,t)$ and $\nabla g : \h \times [0,1]
\rightarrow \h$ is a compact mapping,
\item[(3)] $\nabla H_{\infty}((v,w),0)= \nabla \Phi(v,w)$
\item [(4)] there exists a $G$-equivariant, gradient mapping
$\nabla \varphi_{\infty} : V_{\infty} \rightarrow V_{\infty}$ such that   $$\nabla H_{\infty}((v,w),1)= (\nabla
\varphi_{\infty}(v), A_{\infty}(w))=(\nabla \varphi_{\infty}(v), (L-L_{\infty})_{|W_{\infty}}(w)). $$ \ee \el
\begin{proof}
Recall, that $\Phi_1^{\infty}= (Id-\pi_{\infty}) \circ \nabla \Phi: \h \rightarrow V_{\infty}$ and
$\Phi_2^{\infty}= \pi_{\infty} \circ \nabla \Phi: \h \rightarrow W_{\infty}$.\linebreak Applying Theorem
\ref{uwiklinfty} to $\Phi$, we obtain $\beta_0
> 0$ and a $G$-equivariant $C^1$-mapping \linebreak $w_{\infty} :
V_{\infty} \setminus cl (D_{\beta_0} (V_{\infty},0)) \rightarrow W_{\infty}$ such that $\Phi_2^{\infty}(v,w) = 0$
for
$v \in V_{\infty} \setminus cl(D_{\beta_0}(V_{\infty},0))$ iff $w=w_{\infty}(v)$.\\
Fix $\beta_1 > \beta_0$ such that
 \beq \label{splitA} (Id - \pi_{\infty})(\nabla
\Phi^{-1}(0)) \subset cl(D_{\beta_1}(V_{\infty},0)) \eeq and
\begin{equation}\label{meanvalueeq}
\sup \{\|\nabla^2\eta_{\infty}(v,w)\|; \  (v,w) \in \h \text{ and }\|(v,w)\|_{\h} > \beta_1\} \leq \frac12
\|A_{\infty}^{-1}\|^{-1}.
\end{equation}
Consider the family $H_{\infty} : (V_{\infty} \oplus W_{\infty}) \times [0,1] \rightarrow \h$ defined by formula
\eqref{homotopia} for $p=\infty$, i.e.
\begin{align*} H_{\infty}((v,w),t)&= \frac12 \langle A_{\infty}(w),w \rangle_{\h} +
\frac 12 t (2-t) \langle A_{\infty}(w_{\infty}(v)), w_{\infty}(v) \rangle_{\h} + \\&+t \eta_{\infty}(v,
w_{\infty}(v)) + (1-t) \eta_{\infty} (v, w+ t w_{\infty}(v)).
\end{align*}
We claim  that there exists $\beta_2 > 0$ such that if $(v,w) \in \h \setminus (cl (D_{\beta_1}(V_{\infty},0))
\times cl (D_{\beta_2}(W_{\infty},0)))$ then $\nabla H_{\infty}((v,w),t) \neq 0.$ Notice that $(v,w) \in \h
\setminus (cl (D_{\beta_1}(V_{\infty},0)) \times cl(D_{\beta_2}(W_{\infty},0))$ iff \\
(i) either $\|v\|_{\h} < \beta_1 \text{ and } \|w\|_{\h} > \beta_2$,   \\
(ii) or $\| v\|_{\h} > \beta_1.$

\nt Case (i). For $\beta_1$ there exists K such that if $\|v\|_{\h}< \beta_1$ then $\|w_{\infty}(v)\|_{\h} \leq
K$. Moreover, from (F.4), for fixed $\varepsilon > 0$ there exists $M > 0$ such that if $\|x\|_{\h} > M$ then
$\|\nabla \eta_{\infty}(x)\|_{\h}< \varepsilon \|x\|_{\h}.$ Taking $\varepsilon < \frac12
\|A_{\infty}^{-1}\|^{-1}$ and sufficiently large $\beta_2$ if $\|w\|_{\h}\geq \beta_2$ then
\begin{equation*}
\begin{split} \| \pi_{\infty}( \nabla H_{\infty}((v,w),t))\|_{\h}&=
\|A_{\infty}(w) +(1-t)\pi_{\infty} \nabla
\eta_{\infty}(v,w+tw_{\infty}(v))\|_{\h} \geq\\
&\geq \|A_{\infty}^{-1}\|^{-1} \|w\|_{\h} -(1-t) \|\nabla \eta_{\infty}(v, w+tw_{\infty}(v))\|_{\h}\geq \\&\geq
\|A_{\infty}^{-1}\|^{-1} \|w\|_{\h}-\varepsilon \|(v,w+tw_{\infty}(v)))\|_{\h} > \\ & > \|A_{\infty}^{-1}\|^{-1}
\|w\|_{\h}-\varepsilon \beta_1 - \varepsilon \|w\|_{\h}-\varepsilon K > \\
& > \frac12\|A_{\infty}^{-1}\|^{-1}\|w\|_{\h}-\varepsilon \beta_1 - \varepsilon K > 0.
\end{split}
\end{equation*}

\nt Case(ii). Fix $(v,t) \in (V_{\infty}\setminus cl(D_{\beta_1}(V_{\infty},0))) \times[0,1].$ We claim that
\begin{equation}\label{caseii}
\pi_{\infty}(\nabla H_{\infty}((v,w),t))=0
\end{equation}
 iff $w=0$. Indeed, $w=0$ is a solution of \eqref{caseii}. We proceed to show that it is
 unique.
\begin{equation*}
\begin{split}
&\pi_{\infty}(\nabla H_{\infty}((v,w),t))=0 \Leftrightarrow\\
&A_{\infty}(w)+(1-t)\pi_{\infty} ( \nabla
\eta_{\infty}(v,w+tw_{\infty}(v)))=0 \Leftrightarrow\\
&w=-(1-t)A_{\infty}^{-1}\circ \pi_{\infty} (\nabla
\eta_{\infty}(v,w+tw_{\infty}(v))).
\end{split}
\end{equation*}
Denote the right side of the above equality by $\phi(v,t)(w).$ From the mean value theorem and \eqref{meanvalueeq}
we obtain

$$\|\phi(v,t)(w_1)-\phi(v,t)(w_2)\|_{\h} \leq \|A_{\infty}^{-1}\| \cdot \|\nabla\eta_{\infty}(v,tw_1+w_{\infty}(v))
-\nabla\eta_{\infty}(v,tw_2+w_{\infty}(v)) \|_{\h} \leq $$
$$\leq \|A_{\infty}^{-1}\| \cdot \|w_1-w_2\|_{\h} \cdot \sup
\{\|\nabla^2\eta_{\infty}(u)\|; \  u \in \h \text{ and }\|u\|_{\h} > \beta_1\} \leq \frac12 \|w_1-w_2\|_{\h}.$$
Hence $\phi(v,t) : W_{\infty} \rightarrow W_{\infty}$ is a contraction. Using the Banach fixed point theorem we
conclude that $w=0$ is the unique solution of \eqref{caseii}. If $w=0$ then $(Id-\pi_{\infty})(\nabla
H_{\infty}((v,0),t))=$ $= \Phi_1^{\infty}(v,w(v)) \neq 0$ for $\|v\| \geq \beta_1.$ Now it suffices to take
$\alpha_{\infty}=max \{\beta_1, \beta_2\}.$\\
Lemma \ref{postac} yields (2). To complete the proof notice that
\begin{enumerate}
    \item[\rm (1)] $\nabla H_{\infty}((v,w),0) = \nabla \Phi(v,w) $,
    \item[\rm (2)] $\nabla H_{\infty}((v,w),1) = (\nabla\varphi_{\infty}(v),A_{\infty}
    (w)),$ where $\nabla \varphi_{\infty}(v) = \Phi_1^{\infty} (v,w_{\infty}(v)).$
\end{enumerate}
\end{proof}

\br Notice that from conditions (1) in Lemmas \ref{losp} and \ref{losinfty} for $p=0, \infty$ respectively, we
have $$\nabla H_p(\cdot, t) : (cl(D_{\alpha_p}(V_p,0)) \times cl(D_{\alpha_p}(W_p,0)),
\partial (cl(D_{\alpha_p}(V_p,0))) \times
cl(D_{\alpha_p}(W_p,0)))) \rightarrow$$ $$\rightarrow (\h, \h \setminus \{0\}).$$ Hence both homotopies $H_p$ are
$\Omega$-admissible for $\Omega=D_{\alpha_p}(V_p,0)\times D_{\alpha_p}(W_p,0).$ \er


\numsec
\section{Applications of Splitting Lemmas}
\label{asl}

In this section we compute the indices of isolated  critical points of $\sone$-invariant functionals in terms of
the degree for $\sone$-equivariant gradient maps. Throughout this section we assume that $L=Id$, $G=\sone$ and
$\dim \: \h^{\sone}<\infty$. Let functional $\Phi \in C^2_{\sone}(\h, \mathbb{R})$ be given by \eqref{funkcjonal}.
Suppose that $\displaystyle{\Phi(x)=\frac12 \langle (L- L_p)x, x \rangle_{\h} + \eta_p(x)}$ for $p \in \{0,
\infty\}$, and   that assumptions $(F.1)-(F.6)$ are satisfied. Notice that if $p=0$, functional $\Phi$ satisfies
assumptions of Lemma \ref{losp} and if $p=\infty,$ assumptions of Lemma \ref{losinfty}. Hence, further
consideration we can carry out paralel for $p \in \{0, \infty\}.$ From   Lemmas \ref{losp}, \ref{losinfty} and
Theorem \ref{pinf} we have
$$\dg(\nabla \Phi,D_{\alpha_p}(\mathbb{H},p))=$$
\begin{equation}\label{stopien1}
 =
 \dg((\nabla\varphi_p,A_p),D_{\alpha_p}(V_p,p)\times
D_{\alpha_p}(W_p,p))=
\end{equation}
$$=\dg(\nabla \varphi_p, D_{\alpha_p}(V_p,p))\star \dg(A_p, D_{\alpha_p}(W_p,p))$$

\bl \label{appllem1}Fix $p \in \{0, \infty\}.$ Let $\Phi \in C^{2}_{\sone}(\mathbb{H},\mathbb{R}),$ admits the
representation $\Phi(x) = \frac{1}{2}\langle(Id-L_p) x,x \rangle_{\h} + \eta_p (x),$ and assumptions $(F.1)-(F.6)$
are fulfilled. Moreover, assume that $V_p=\ker(Id-L_p)\subset \h^{\sone}.$ Then
$$\nabla_{\sone}-\mathrm{deg}(\nabla\Phi,D_{\alpha_p}(\h,p)) =  \mathrm{deg}_{\rm
B}(\nabla\Phi^{\sone},D_{\alpha_p}(\h^{\sone},p),0) \cdot$$
$$\cdot \dg((Id-L_p)_{\mid (\h^{\sone})^{\bot} },D_{\alpha_p}((\h^{\sone})^{\bot},0)).$$
 \el
\begin{proof} Since $A_p=((A_p)_{|W_p^{\sone}},
(A_p)_{|(W_p \ominus W_p^{\sone})}),$ $V_p \subset \h^{\sone}$ and \eqref{stopien1} we have
$$\dg(\nabla \Phi, D_{\alpha_p}(\h,p))=$$ $$=\dg (\nabla \varphi_p,
D_{\alpha_p}(V_p,p)) \star   \dg ((A_p)_{|W_p^{\sone}}, D_{\alpha_p}(W_p^{\sone},0)) \star$$
$$\star \dg ((A_p)_{|(W_p \ominus W_p^{\sone})}, D_{\alpha_p}((W_p \ominus W_p^{\sone}),0))=$$

$$ =\dg(( \nabla \varphi_p,(A_p)_{|(W_p \ominus W_p^{\sone})}),D_{\alpha_p}(V_p,p)\times D_{\alpha_p}((W_p \ominus W_p^{\sone}),0) \star $$

$$\star \dg ((A_p)_{|(\h^{\sone})^{\bot}}, D_{\alpha_p}((\h^{\sone})^{\bot},0))= \dg (\nabla \Phi^{\sone},
D_{\alpha_p}(\h^{\sone},p))\star $$
\begin{equation}\label{opa}
 \star \dg((Id-L_p)_{\mid (\h^{\sone})^{\bot} },D_{\alpha_p}((\h^{\sone})^{\bot},0))
\end{equation}

\nt  Recall that  $\big(\nabla \Phi\big)^{\sone}=\big(\nabla \Phi\big)_{|\h^{\sone}}=\nabla \big(
\Phi_{|\h^{\sone}}\big) : \h^{\sone} \rightarrow \h^{\sone}$ is well-de\-fi\-ned gradient map and notice that by
Remark \ref{goodrem} we have
\begin{equation}\label{opass}
  \dg(\nabla \Phi^{\sone},D_{\alpha_p}(\h^{\sone},0))= (\mathrm
{deg}_{\rm B}(\nabla\Phi^{\sone},D_{\alpha_p}(\h^{\sone},0)), 0,\dotsc )
\end{equation}
Taking into account formulas \eqref{m}, \eqref{opa} and \eqref{opass} we complete the proof.
\end{proof}

\noindent Recall, that by $V_L(\lambda_i)$ we denote the eigenspace of the operator L corresponding to the
eigenvalue $\lambda_i$. Using Theorem \ref{dizom} we obtain
\begin{equation*}
\begin{split}
&\dg(\nabla \Phi,D_{\alpha_p}(\h,p)) =
\\
&=\dg (\nabla \varphi_p, D_{\alpha_p}(V_p,p)) \star \dg(-Id, D_{\alpha_p} (\bigoplus_{\lambda_i
> 1}V_{L_p}(\lambda_i),0)).
\end{split}
\end{equation*}
\noindent Since $V_p = \ker (Id-L_p)$ is a finite-dimensional representation of the group $\sone$, from Theorem
\ref{tk}  we obtain    numbers $j_0 \geq 0,  j_1, \dotsc, j_r > 0$ and $k_r > \dotsc > k_1 > k_0=0 $ such that
$V_p \approx \bR[j_0,k_0] \oplus \mathbb{R}[j_1,k_1] \oplus \dotsc \oplus \mathbb{R}[j_r,k_r]$. Let $\{i_1,
\dotsc, i_s \} \subset \{1, \dotsc ,r\}$. Denote   $k_{i_1 \dotsc i_s}=\gcd\{k_{i_1}, \dotsc, k_{i_s}\}$. Possible
isotropy groups of points of $V_p$ are $\sone$ and groups $\mathbb{Z}_{k_{i_1 \dotsc i_s}}$, for arbitrary $\{i_1,
\dotsc, i_s\} \subset \{1, \dotsc, r\}.$ From Remark \ref{goodrem}  we have $\dg_H(\nabla \varphi_p,
D_{\alpha_p}(V_p,p))=0$ for $H \notin \{\sone\} \cup \displaystyle{ \bigcup_{\{i_1, \dotsc, i_s\} \subset \{1,
\dotsc ,r\}}} \{\mathbb{Z}_{k_{i_1 \dotsc i_s}}\}$. Hence we have proved the following lemma.

\bl\label{ogolny} Let $p \in \{0, \infty\}.$ Let $\Phi \in C^{2}_{\sone}(\mathbb{H},\mathbb{R}),$ admits the
representation $\Phi(x) = \frac{1}{2}\langle(Id-L_p) x,x \rangle_{\h} + \eta_p (x),$ and assumptions $(F.1)-(F.6)$
hold. Then
\begin{equation*}
\begin{split}
&\dg(\nabla \Phi,D_{\alpha_p}(\h,p)) =
\\
&=\dg (\nabla \varphi_p, D_{\alpha_p}(V_p,p)) \star \dg(-Id, D_{\alpha_p}(\bigoplus_{\lambda_i > 1}
V_{L_p}(\lambda_i),0)).
\end{split}
\end{equation*}
Moreover, by formula \eqref{m}, we have
\begin{equation}\label{ogolnydeg}
\begin{split}
& \displaystyle \dg_H(\nabla \Phi, D_{\alpha_p}(\h,p))=\\&= \dg_{\sone}(\nabla \varphi_p,
D_{\alpha_p}(V_p,p))\cdot \dg_H(-Id, D_{\alpha_p}(\bigoplus_{\lambda_i > 1} V_{L_p}(\lambda_i),0))
\end{split}
\end{equation}
for $H \notin \{\sone \} \cup \displaystyle{ \bigcup_{\{i_1, \dotsc, i_s\} \subset \{1, \dotsc ,r\}}}
\{\mathbb{Z}_{k_{i_1 \dotsc i_s}}\}.$ \el

\begin{Corollary}\label{wnappl}
Let assumptions of Lemma \ref{ogolny} be satisfied. Assume additionally that
$V_p^{\sone}=\ker(Id-L_p)^{\sone}=\{0\}.$ Then
\begin{equation*}
\ds \dg_H (\nabla \Phi, D_{\alpha_p}(\h,p))= \dg_H(-Id, D_{\alpha_p}(\bigoplus_{\lambda_i > 1}
V_{L_p}(\lambda_i),0))
\end{equation*}
for $H \notin \{\sone\} \cup \displaystyle{ \bigcup_{\{i_1, \dotsc, i_s\} \subset \{1, \dotsc ,r\}}}
\{\mathbb{Z}_{k_{i_1 \dotsc i_s}}\}.$
\end{Corollary}

\begin{proof}
Since $V_p^{\sone}=\ker(Id-L_p)^{\sone}=\{0\}$, by Remark \ref{goodrem} we have
$$\dg_{\sone}(\nabla \varphi_p, D_{\alpha_p}(V_p,p))=1.$$
Applying formula \eqref{ogolnydeg} we complete the proof.
\end{proof}


\section{Nonstationary Periodic Solutions of Autonomous Newtonian Systems}
\label{persol}

Throughout this section we study periodic solutions of autonomous Newtonian systems. We define an
$\sone$-invariant functional  on a suitably chosen  infinite-dimensional Hilbert space which is an
infinite-dimensional, orthogonal representation  of the group $\sone.$ Critical orbits of this functional are in
one-to-one correspondence with solutions of a considered system. Therefore for our purpose it is enough to study
only the critical orbits of this functional.

\nt We begin this section with a definition of an appropriate Hilbert space. Fix $T > 0$ and define
$$\h^1_T = \{u : [0,T] \rightarrow \bR^n : \text{ u is abs. cont., } u(0)=u(T), \dot u \in
L^2([0,T],\bR^n)\}.$$ It is known that $\h^1_T$ is a separable Hilbert space with a scalar product given by the
formula $\ds \langle u,v\rangle_{\h^1_T} = \int_0^T (\dot u(t), \dot v(t)) + (u(t),v(t)) \; dt,$ where $(\cdot,
\cdot)$ and $\| \cdot \|$ are the usual scalar product and norm in  $\bR^n,$ respectively. It is easy to show that
$\left(\h^1_T,\langle\cdot,\cdot\rangle_{\h^1_T}\right)$ is an orthogonal representation of the group $\sone$ with
an $\sone$-action given by shift in time.

\nt Let us consider the following Newtonian system

\begin{equation}\label{NNS}
 \begin{cases}
    \ddot u =  -V'(u)&  \\
    u(0)=u(T) & \\
    \dot u (0)=\dot u(T)
 \end{cases}
\end{equation}
where $V \in C^2(\bR^n,\bR).$  Solutions of  \eqref{NNS} are in one to one correspondence with critical points of
an $\sone$-invariant $C^2$-functional  $\Phi_V : \h^1_{T} \rightarrow \bR$ defined as follows

\begin{equation}\label{NFUN}
  \Phi_V(u) = \int_0^T \frac{1}{2} \| \dot u(t) \|^2 - V(u(t)) \; dt
\end{equation}

\nt Notice that for any $u,v \in \h^1_T$ we have $\Phi_V'(u)(v)=\langle \nabla\Phi_V (u),v\rangle_{\h^1_T}=
\langle u - \nabla \zeta(u),v\rangle_{\h^1_T},$ where $\nabla \zeta : \h^1_T \rightarrow \h^1_T$ is an
$\sone$-equivariant, compact, gradient operator given by the formula  $\ds \langle \nabla
\zeta(u),v\rangle_{\h^1_T}=\int_0^T (u(t)+V(u(t)),v(t)) \; dt.$ In the other words the gradient $\nabla \Phi_V :
\h^1_T \rightarrow \h^1_T$ is an $\sone$-equivariant $C^1$-operator in the form of a compact perturbation of the
identity.

\numsubsec
\subsection{Linear Equation}
\label{persollin}

 In this section we carry on  detailed analysis of a linear system of the form
\begin{equation}\label{LNS}
 \begin{cases}
    \ddot u =  -A u&  \\
    u(0)=u(T) & \\
    \dot u (0)=\dot u(T)
 \end{cases}
\end{equation}
where $A$ is a real, symmetric $(n \times n)$-matrix. Moreover, we study properties of a functional associated
with equation \eqref{LNS}.

\nt Define the corresponding functional  $\Phi_A : \h^1_{T} \rightarrow \bR$  as follows

\begin{equation}\label{LFUN}
  \Phi_A(u) = \frac{1}{2} \int_0^T  \| \dot u(t) \|^2 - (Au(t),u(t)) \; dt =\frac{1}{2} \langle u-L_A(u),u \rangle_{\h^1_T}
\end{equation}
and notice that $\nabla \Phi_A = Id - L_A,$ where $L_A : \h^1_T \rightarrow \h^1_T$ is a linear, self-adjoint,
$\sone$-e\-qui\-va\-riant and compact operator defined by the formula $$\langle L_A (u),v\rangle_{\h^1_T}=\ds
\int_0^T (u(t)+Au(t),v(t)) dt.$$

\nt To study the linear eigenvalue problem $u - \lambda L_{Id_{\bR^n}} (u)=0$ it is enough to consider the
following system
\begin{equation}\label{ELNS}
 \begin{cases}
    \ddot u =  -(2 \lambda-1) u&  \\
    u(0)=u(T) & \\
    \dot u (0)=\dot u(T)
 \end{cases}
\end{equation}

\nt It is easy to check that eigenvalues and eigenspaces of the operator $L_{Id_{\bR^n}} : \h^1_T \rightarrow
\h^1_T$ are of the following form
\begin{enumerate}
  \item $\ds \sigma\left(L_{Id_{\bR^n}}\right)=\left\{\lambda_k=\frac{2T^2}{T^2+4 k^2 \pi^2}\right\}_{k \in \bN \cup \{0\}},$
  \item $\ds V_{L_{Id_{\bR^n}}}(\lambda_0)=\bR^n \approx \bR[n,0],$
  \item $\ds V_{L_{Id_{\bR^n}}}\left(\lambda_k\right)=
  \left\{a_k \cos \frac{2 k \pi}{T} t + b_k \sin \frac{2 k \pi}{T} t : a_k,b_k \in \bR^n\right\} \approx \bR[n,k],$
 $k \in \bN.$
\end{enumerate}
From the above we obtain that an orthonormal basis in $\h^1_T$ can be chosen as follows

$$\ds\sqrt{\frac{1}{T}} \cdot e_i, \sqrt{\frac{2T}{T^2+4k^2 \pi^2}} \cdot \cos \left(\frac{2 k \pi}{T} t\right)
\cdot e_i,\sqrt{\frac{2T}{T^2+ 4k^2 \pi^2}} \cdot \sin \left(\frac{2 k \pi}{T} t\right)  \cdot e_i$$ or in
equivalent way
$$\ds \sqrt{\frac{\lambda_0}{2T}} \cdot e_i,  \sqrt{\frac{\lambda_k}{T}} \cdot \cos \left(\frac{2 k \pi}{T} t\right)
\cdot e_i, \sqrt{\frac{\lambda_k}{T}} \cdot \sin \left(\frac{2 k \pi}{T} t\right)  \cdot e_i,$$ where
$i=1,\ldots,n$ and $k \in \bN.$

\nt It is clear that $u \in \h^1_T$ possesses  Fourier series of the form
$$u(t)= \tilde{a}_0 \sqrt{\frac{\lambda_0}{2T}} + \sum_{k \in \bN} \tilde{a}_k \cdot \bigg(\sqrt{\frac{\lambda_k}{T}} \cdot \cos \frac{2 k \pi}{T} t \bigg) +
 \tilde{b}_k \cdot \bigg(\sqrt{\frac{\lambda_k}{T}} \cdot \sin \frac{2 k \pi}{T} t\bigg)=$$
 $$=a_0 + \sum_{k \in \bN} a_k \cdot \cos \frac{2 k \pi}{T} t +
 b_k  \cdot \sin \frac{2 k \pi}{T} t.$$

\nt In the following lemma we study properties of the operator $\nabla \Phi_A = Id - L_A : \h^1_T \rightarrow
\h^1_T.$

\bl \label{lgrad} If $u \in \h^1_T$ with Fourier series
$$u(t)= a_0   + \sum_{k \in \bN} a_k \cdot   \cos \frac{2 k \pi}{T} t   +
 b_k \cdot   \sin \frac{2 k \pi}{T} t,$$
then

$\ds\nabla \Phi_A(u)= u-L_A (u)=-A \cdot a_0  + \sum_{k=1}^{\infty}
 (\Lambda(k) \cdot a_k)   \cdot \cos \frac{2 k \pi}{T} t +  (\Lambda(k) \cdot b_k)  \cdot \sin \frac{2 k
\pi}{T} t,$

\nt where $\ds \Lambda(k)=\left(\frac{4 k^2 \pi^2}{4 k^2 \pi^2 + T^2} Id - \frac{T^2}{4 k^2 \pi^2 + T^2} A\right)$
\el
\begin{proof} Fix $u \in \h^1_T$ with Fourier series $\ds u(t)= a_0 + \sum_{k \in \bN} a_k \cdot \cos \frac{2 k \pi}{T} t +
 b_k  \cdot \sin \frac{2 k \pi}{T} t.$ Then
$$\nabla \Phi_A(u)=u-L_A (u)= (Id -L_A)( a_0)+  \sum_{k \in \bN} (Id-L_A)\left(a_k \cdot \cos \frac{2 k \pi}{T} t +
b_k \cdot \sin \frac{2 k \pi}{T} t  \right).$$ What is left is to
compute
\begin{enumerate}
  \item $\ds (Id-L_A)( a_0) \in V_{L_{Id_{\bR^n}}}(\lambda_0),$
  \item $\ds (Id-L_A)\left(a_k \cdot \cos \frac{2 k \pi}{T} t + b_k \cdot \sin \frac{2 k
\pi}{T} t \right) \in   V_{L_{Id_{\bR^n}}}(\lambda_k).$
\end{enumerate}
Put $u_0=a_0,$   fix $v \in \h^1_T$ and notice that
$$\langle\nabla \Phi_A(u_0),v \rangle_{\h^1_T}=\langle u_0 - L_A(u_0),v \rangle_{\h^1_T}= \int_0^T  - (A (u_0) ,v(t))\;dt=$$
$$= \int_0^T  (\dot u_0 , \dot v(t)) +(u_0,v(t))- ((Id+A) (u_0) ,v(t))\;dt =$$
$$=\langle u_0 - (Id + A)(u_0),v\rangle_{\h^1_T}=\langle  -   A(u_0),v\rangle_{\h^1_T}.$$

\nt Summing up, we obtain $\ds \nabla \Phi_A(u_0)=u_0-L_A(u_0)=-A
(u_0).$

\nt For simplicity of notation, we let  $\ds u_k(t)$ stand for
$\ds a_k \cdot \cos \frac{2 k \pi}{T} t + b_k \cdot \sin \frac{2 k
\pi}{T} t.$ Fix $v \in \h^1_T$ and notice that

$$\langle\nabla \Phi_A(u_k),v \rangle_{\h^1_T}=\langle u_k - L_A(u_k),v \rangle_{\h^1_T}=
\int_0^T \; (\dot u_k(t), \dot v(t))- (A u_k(t),v(t))\;dt=$$

$$= \int_0^T (\dot u_k(t), \dot v(t)) + (u_k(t),v(t))- (u_k(t)+A u_k(t),v(t))\; dt=
\int_0^T (-\ddot u_k(t)+u_k(t),v(t))- $$ $$-(u_k(t)+A u_k(t),v(t))\; dt= \int_0^T \left(\frac{4 k^2
\pi^2+T^2}{T^2}\right)(u_k(t),v(t)) \; dt-$$ $$- \int_0^T  \left(\frac{T^2}{4 k^2
\pi^2+T^2}\right)\left(\left(\frac{4 k^2 \pi^2+T^2}{T^2}\right) \left(u_k(t)+A u_k(t)\right),v(t)\right)\; dt=$$
$$=\langle u_k,v \rangle_{\h^1_T} - \left\langle \left(\frac{T^2}{4 k^2
\pi^2+T^2}\right)(Id+A) u_k,v\right\rangle_{\h^1_T}=$$
$$=\left\langle \left(\frac{4 k^2 \pi^2}{4 k^2 \pi^2 + T^2}
Id - \frac{T^2}{4 k^2 \pi^2 + T^2}
A\right)u_k,v\right\rangle_{\h^1_T}.$$

\nt Summing up, we obtain
$$\nabla \Phi_A(u_k)=u_k-L_A(u_k)=\left(\frac{4 k^2 \pi^2}{4 k^2 \pi^2 + T^2}
Id - \frac{T^2}{4 k^2 \pi^2 + T^2} A\right)u_k=$$
$$= (\Lambda(k) \cdot a_k)   \cdot \cos \frac{2 k \pi}{T} t +  (\Lambda(k) \cdot b_k)  \cdot \sin \frac{2 k
\pi}{T} t= \Lambda(k) u_k,$$ which completes the proof.
\end{proof}

\nt As a direct consequence of Lemma \ref{lgrad} we obtain the
following two corollaries.

\bco \label{wniosek1} The following conditions are equivalent
\begin{enumerate}
  \item operator $\nabla \Phi_A=Id - L_A : \h^1_T \rightarrow \h^1_T$ is an isomorphism,
  \item $\ker \Lambda(k)=\{0\} $ for any $k \in \bN \cup \{0\},$
  \item $\ds \sigma(A) \cap
\left\{\frac{4 k^2 \pi^2}{T^2} : k \in \bN \cup \{0\}\right\}=
\emptyset.$
\end{enumerate}
\eco

\nt For $\alpha \in \bR$  we will denote by $\mu_A(\alpha)$ the multiplicity of $\alpha$ considered as an
eigenvalue of matrix $A.$ If $\alpha \notin \sigma(A) $ then it is understood that $\mu_A(\alpha)=0.$ Moreover, if
$\alpha \in \sigma(A)$ then  we will denote by $V_A(\alpha)$ the eigenspace of $A$ corresponding to the eigenvalue
$\alpha.$ For any  $k \in \bN \cup \{0\}$ define
\begin{enumerate}
  \item $\ds \sigma_k(A,T)=\sigma(A) \cap\left(\frac{4k^2\pi^2}{T^2},+\infty \right),$
  \item $\ds j_k(A,T)=\sum_{\alpha \in \sigma_k(A,T)} \mu_A(\alpha).$
\end{enumerate}

\br  Notice that
$$\nu(A,T)= \mu_A(0) +2\sum_{k=1}^{\infty} \mu_A\left( \frac{4k^2\pi^2}{T^2}\right),
j(A,T)=j_0(A,T)+2 \sum_{k=1}^{\infty} j_k(A,T),$$ where the numbers $\nu(A,T), j(A,T)$ are defined in
\cite{[MAW-WIL]} on page 207. Since we are going to apply the degree for $\sone$-equivariant gradient maps, we
have to describe the finite-dimensional spaces $\ker Id - L_A$ and  $\ds \bigoplus_{\lambda_i > 1}
V_{L_A}(\lambda_i)$ as representations of the group $\sone.$ It is not enough for our purpose  to know only the
dimensions $\nu(A,T)= \dim \ker Id - L_A$ and $\ds j(A,T) = \dim  \bigoplus_{\lambda_i > 1}
V_{L_A}(\lambda_i).$\er

\nt The following corollary will prove extremely useful in the next sections.

\bco\label{wniosek} Operator $\nabla \Phi_A=Id - L_A : \h^1_T \rightarrow \h^1_T$ has the following properties
\begin{enumerate}
  \item $\ds \h^1_{T,0}=\ker \nabla \Phi_A = V_{L_A}(1)=$
  $$ = \ker A \oplus \bigoplus_{k=1}^{\infty}     \left\{a_k \cdot \cos \frac{2 k \pi}{T} t + b_k \cdot \sin \frac{2 k
\pi}{T} t : a_k,b_k \in V_A\left(\frac{4k^2\pi^2}{T^2}\right)\right\},$$
  \item $\ds \h^1_{T,0} \approx   \bigoplus_{k=0}^{\infty}
  \bR\left[\mu_A\left( \frac{4k^2\pi^2}{T^2}\right),k\right],$
  \item $\ds \dim \h^1_{T,0}=  \mu_A(0) +2\sum_{k=1}^{\infty} \mu_A\left( \frac{4k^2\pi^2}{T^2}\right),$
  \item $\ds \h^1_{T,-}= \bigoplus_{\lambda_i > 1} V_{L_A}(\lambda_i)=$
  $$=\bigoplus_{\alpha \in \sigma_0(A)} V_A(\alpha) \oplus
  \bigoplus_{k=1}^{\infty} \bigoplus_{\alpha \in \sigma_k(A,T)}  \left\{a_k \cdot \cos \frac{2 k \pi}{T} t + b_k \cdot \sin \frac{2 k
\pi}{T} t : a_k,b_k \in V_A\left(\alpha\right)\right\},$$
  \item $\ds \h^1_{T,-} \approx \bigoplus_{k=0}^{\infty} \bigoplus_{\alpha \in \sigma_k(A,T)}
  \bR\left[\mu_A\left( \alpha \right),k\right]=\bigoplus_{k=0}^{\infty} \bR[j_k(A,T),k],$
  \item $\ds \dim \h^1_{T,-}=\sum_{\alpha \in \sigma_0(A)} \mu_A(\alpha) +
  2 \sum_{k=1}^{\infty} \sum_{\alpha \in \sigma_k(A,T)} \mu_A(\alpha)=j_0(A,T)+2 \sum_{k=1}^{\infty} j_k(A,T),$
  \item $\ds \h^1_{T,+}= \overline{\bigoplus_{\lambda_i <1} V_{L_A}(\lambda_i)},$
  \item $\h^1_T= \h^1_{T,-} \oplus \h^1_{T,0} \oplus \h^1_{T,+}.$
\end{enumerate}
\eco

\nt The following fact is a  direct consequence of Theorem \ref{dizom}, Lemma \ref{lgrad} and Corollaries
\ref{wniosek1}, \ref{wniosek}.

\bft\label{faktdegizo} If $\ds \sigma(A) \cap \left\{\frac{4 k^2 \pi^2}{T^2} : k \in \bN \cup \{0\}\right\}=
\emptyset,$ then the operator $\nabla \Phi_A= Id - L_A : \h^1_T \rightarrow \h^1_T$ is an isomorphism.
 Additionally, for $\gamma > 0$
$$\nabla_{\sone}-\mathrm{deg}_H\left(\nabla \Phi_A,D_{\gamma}\left(\h^1_T,0\right)\right)=
\nabla_{\sone}-\mathrm{deg}_H\left(-Id,D_{\gamma}\left(\h^1_{T,-},0\right)\right)=$$
$$= \left\{\begin{array}{ll}
  (-1)^{j_0(A,T)}  &    \text{ for } H=\sone,\\
  & \\
  (-1)^{j_0(A,T)} \cdot j_k(A,T) &  \text{ for }  H=\bZ_{ k}.
\end{array} \right.$$
It is understood that if $\ds \h^1_{T,-}=\{0\}$ then $\ds \nabla_{\sone}-\mathrm{deg}\left(\nabla
\Phi_A,D_{\gamma}\left(\h^1_T,0\right)\right)= \bI \in U(\sone).$\eft

\numsubsec
\subsection{Existence of Periodic Solutions of Nonlinear Equation}
\label{persolnonlin}

In this section we formulate sufficient conditions for the existence of non-stationary $T$-periodic solutions of
autonomous Newtonian systems.

 \nt Let us consider the following Newtonian system

\begin{equation}\label{NNSO}
 \begin{cases}
    \ddot u =  -V'(u)&  \\
    u(0)=u(T) & \\
    \dot u (0)=\dot u(T)
 \end{cases}
\end{equation}
where $V \in C^2(\bR^n,\bR).$ Suppose
\begin{enumerate}
\item[(i)] $(V')^{-1}(0)= \{p_1, \dotsc, p_q\}$,
\item[(ii)] $V'(x)= V''(\infty) \cdot x+o(\|x\|)$ as $\|x\| \rightarrow \infty.$
\end{enumerate}
Recall that for $p \in \h^{\sone}$ and $\gamma > 0$ we set $D_{\gamma}(\h,p)=\{v \in \h : \|v-p\|_{\h} <
\gamma\}.$ Moreover, if $p=\infty$ then $D_{\gamma}(\h,\infty):=D_{\gamma}(\h,0).$ Define $\ds
\ind(-V',p_i)=\lim_{\alpha \rightarrow 0 }\deg_{\rm B}(-V',D_{\alpha}(\bR^n,p_i),0)$ for $i=1,\ldots,q$ and $\ds
\ind(-V',\infty)=\lim_{\alpha \rightarrow \infty}\deg_{\rm B}(-V',D_{\alpha}(\bR^n,\infty),0).$

\br \label{sumfo} Under the above assumptions  $\ds \ind(-V',\infty)=\sum_{i=1}^q \ind(-V',p_i).$ \er

\begin{Definition}\label{Defindex}
Let $p\in \{p_1,\ldots,p_q,\infty\}$ and $\ds \sigma(V''(p)) \cap \left\{ \frac{4k^2\pi^2}{T^2} : \ k \in \bN
\right\}= \emptyset$. Define $I_V(p,T)=(I_V(p,T)_{\sone},I_V(p,T)_{\bZ_1}, \dotsc, I_V(p,T)_{\bZ_k}, \dotsc,) \in
U(\sone)$ in the following way

\begin{equation}\label{indekscase1}
I_V(p,T)_H=
 \left\{\begin{array}{ll}
  \mathrm{ind}(-V',p))   &    \text{ for } H=\sone,\\
  & \\
  \mathrm{ind}(-V',p)) \cdot j_k\left(V''(p),T\right) &  \text{ for }  H=\bZ_k.
\end{array} \right.
\end{equation}
\end{Definition}

\begin{Remark} Notice that, if moreover we assume  $\det V''(p)  \neq 0,$ then formula \eqref{indekscase1} becomes
\begin{equation}\label{indekscase3}
I_V(p,T)_H=
 \left\{\begin{array}{ll}
  (-1)^{j_0(V''(p),T)}  &    \text{ for } H=\sone,\\
  & \\
  (-1)^{j_0(V''(p),T)}  \cdot j_k\left(V''(p),T\right) &  \text{ for }  H=\bZ_k.
\end{array} \right.
\end{equation}

\end{Remark}

\begin{Lemma}\label{redukcja}
If $p\in \{p_1,\ldots,p_q,\infty\}$ and $\ds \sigma(V''(p)) \cap \left\{ \frac{4k^2\pi^2}{T^2} : \ k \in \bN
\right\}= \emptyset$ then the stationary solution $p$ is an isolated critical point of $\Phi_V$.
\end{Lemma}
\begin{proof}
Since $\ds \sigma(V''(p)) \cap \left\{ \frac{4k^2\pi^2}{T^2} : \ k \in \bN \right\}= \emptyset,$ from Corollary
\ref{wniosek} we conclude that $\ker(Id-L_p)\subset (\h^1_T)^{\sone}.$ Suppose, contrary to our claim that $p$ is
not isolated in $(\nabla \Phi)^{-1}(0).$ Taking into account that $\ker(Id-L_p)\subset (\h^1_T)^{\sone}$ and
Corollary \ref{cpi} we obtain that  $p$ is not isolated in $(\nabla \Phi)^{-1}(0) \cap (\h^1_T)^{\sone}.$   But
there is only a finite number of stationary solutions of \eqref{NNSO}, a contradiction.
\end{proof}

\begin{Lemma}\label{I=deg}
If $p\in \{p_1,\ldots,p_q,\infty\}$ and $\ds \sigma(V''(p)) \cap \left\{ \frac{4k^2\pi^2}{T^2} : \ k \in \bN
\right\}= \emptyset$ then $$I_V(p,T)= \dg(\nabla \Phi_V, D_{\alpha_p}(\h^1_T,p)),$$ where $\alpha_p$ is
sufficiently small if $p=p_i$ for  $i=1,\ldots,q$ and sufficiently large if $p=\infty.$
\end{Lemma}
\begin{proof}
From Lemma \ref{redukcja} it follows that the functional $\Phi_V$ satisfies conditions (F.1)-(F.6). Hence we can
apply Lemma \ref{appllem1} to compute $\dg(\nabla \Phi_V,D_{\alpha_p}(\h^1_T,p)).$ We obtain

$$\dg(\nabla \Phi_V,D_{\alpha_p}(\h^1_T,p)) =
\mathrm{deg}_B(\nabla\Phi_V^{\sone},D_{\alpha_p}((\h^1_T)^{\sone},p),0)\cdot$$
$$\cdot \dg((Id-L_p)_{\mid ((\h^1_T)^{\sone})^{\bot} },D_{\alpha_p}(((\h^1_T)^{\sone})^{\bot},0)).$$
Since  $(Id- L_p)_{\mid ((\h^1_T)^{\sone})^{\bot} }$ is an isomorphism, from Fact \ref{faktdegizo} we have
\begin{equation*}
\begin{split}
\dg((Id -L_p)_{\mid ((\h^1_T)^{\sone})^{\bot} },D_{\alpha_p}(((\h^1_T)^{\sone})^{\bot},0))=\\=(1,
j_1(V''(p),T),\dotsc, j_k(V''(p),T), \dotsc).&
\end{split}
\end{equation*}
Moreover, since $(\h^1_T)^{\sone}=\{ u \in \h^1_T : \ u=const\}$ it is easily seen that
$$\mathrm{deg}_B\big(\nabla\Phi_V^{\sone},D_{\alpha_p}\big(\big(\h^1_T\big)^{\sone},p\big),0\big)
=\mathrm{deg}_B(-V', D_{\alpha_p}(\mathbb{R}^n,p), 0).$$ Applying \eqref{m} we obtain the assertion.
\end{proof}

\nt The following  theorem ensures the existence of at least one non-stationary $T$-periodic solution of system
\eqref{NNSO}. See  \cite{[RADZ],[RARY]} for related results.

\begin{Theorem}\label{istnienie}
Suppose that
\begin{enumerate}
\item[(i)] $(V')^{-1}(0)= \{p_1, \dotsc, p_q\}$,
\item[(ii)] $V'(x)= V''(\infty) \cdot x+o(\|x\|)$ as $\|x\| \rightarrow \infty.$
\end{enumerate}
Assume additionally that  $\ds \sigma(V''(p)) \cap \left\{ \frac{4k^2\pi^2}{T^2} : \ k \in \bN \right\}=
\emptyset$ for any $p\in \{p_1,\dotsc,p_q, \infty\}$ and that $\ds I_V(\infty,T)\neq  \sum_{i=1}^qI_V(p_i,T).$
Then there exists at least one non-stationary $T$-pe\-rio\-dic solution of \eqref{NNSO}.
\end{Theorem}
\begin{proof} Notice that stationary solutions $p_1,\dotsc,p_q$ and
$\infty$ are isolated critical points of functional $\Phi_V$ given by \eqref{NFUN}. Suppose, contrary to our
claim, that $p_1,  \dotsc, p_q$ are the only $T$-periodic solutions of \eqref{NNSO}. Thus we can choose
$\alpha_{\infty},\alpha_{p_i} > 0, i=1,\dotsc,q$ such that
\begin{enumerate}
\item[(i)] $\big(\nabla \Phi_V\big)^{-1}(0) \cap(\h^1_T \setminus
D_{\alpha_{\infty} }(\h^1_T,\infty))= \emptyset,$
\item[(ii)] $\big(\nabla \Phi_V \big)^{-1}(0) \cap D_{\alpha_{p_i}}(\h^1_T,p_i)
= \{p_i\}, i= 1,\dotsc,q$,
\item[(iii)] $D_{\alpha_{p_i}}(\h^1_T,p_i) \cap D_{\alpha_{p_j}}(\h^1_T,p_j)=
\emptyset$, for $i \neq j,$
\item[(iv)] $cl(D_{\alpha_{p_i}}(\h^1_T,p_i)) \subset
D_{\alpha_{\infty}}(\h^1_T,\infty)$ for $i= 1, \dotsc , q.$
\end{enumerate}
From the Theorem \ref{wlas}, \eqref{w3} and \eqref{w4} we obtain
$$\dg(\nabla \Phi_V, D_{\alpha_{\infty}}(\h^1_T,\infty)) =
\sum_{i=1}^q \dg (\nabla \Phi_V, D_{\alpha_{p_i}}(\h^1_T,p_i)).$$ But since $\Phi_V$ satisfies assumptions of the
Lemma \ref{I=deg}, we have
$$\dg(\nabla \Phi_V, D_{\alpha_{\infty}}(\h^1_T,\infty))= I_V(\infty,T) \neq$$
$$\neq \sum_{i=1}^qI_V(p_i,T)= \sum_{i=1}^q \dg (\nabla \Phi_V, D_{\alpha_{p_i}}(\h^1_T,p_i)),$$
a contradiction.
\end{proof}

\begin{Remark}\label{wnistnienie1}
Notice that the above theorem  can be formulated in the following  e\-qui\-va\-lent way. Assume that for $\ds p
\in \{p_1, \dotsc, p_q,\infty\}, \ \sigma(V''(p)) \cap \left\{ \frac{4k^2\pi^2}{T^2} : \ k \in \bN \right\}=
\emptyset$. If there exists $k \in \bN$ such that
$$\ind(-V',\infty)  \cdot
j_k\left(V''(\infty),T\right) \neq \displaystyle{\sum_{i=1}^q \ind(-V',p_i) \cdot j_k\left(V''(p_i),T\right)}$$
then there exists at least one non-stationary $T$-periodic solution of \eqref{NNSO}.
\end{Remark}

\nt Suppose now that $\ds \sigma (V''(p)) \cap \left\{\frac{4k^2\pi^2}{T^2} : k \in \bN\right\} \neq \emptyset.$
In that case we cannot define $I_V(p,T)$ such as in Definition \ref{Defindex}. However, we can define almost all
the coordinates of it. Recall that for $\{k_1,\dotsc,k_r\}$ and $ \{i_1,\dotsc,i_s\} \subset \{1,\dotsc,r\}$ we
have $k_{i_1 \dotsc i_s}=\gcd\{k_{i_1},\dotsc,k_{i_s}\}.$

\begin{Definition}\label{indeksdef2}
Let $\ds  \sigma(V''(p)) \cap \left\{ \frac{4k^2\pi^2}{T^2} : \ k \in \bN \right\} = \left\{ \frac{4
k_1^2\pi^2}{T^2}, \dotsc,\frac{4 k_r^2\pi^2}{T^2}\right\},$ where  $p\in \{p_1,\ldots,p_q, \infty\}$. Put $\ds
\bK = \bigcup_{\{i_1, \dotsc, i_s\} \in \{1, \dotsc, r\}} \{k_{i_1\dotsc i_s}\}$ and define
\begin{equation}\label{indekscase4} I_V(p,T)_H=
 \left\{\begin{array}{ll}
  \ind(-V',p)  &    \text{ for } H=\sone,\\
  & \\
   \ind(-V',p)  \cdot j_k(V''(p),T)  &  \text{ for }  H=\bZ_k \text{ and } k \notin \bK.
\end{array} \right.
\end{equation}
\end{Definition}

\begin{Remark}
Observe that if $\ds \sigma(V''(p)) \cap \left\{ \frac{4k^2\pi^2}{T^2} : \ k \in \bN  \right\} = \left\{
\frac{4k_1^2\pi^2}{T^2}, \dotsc,\frac{4k_r^2\pi^2}{T^2}\right\}$ and moreover $\det V''(p)\neq0$, then
$\ds\deg_B(-V',D_{\alpha_p}(\bR^n,p),0)=(-1)^{j_0(V''(p))}$ and consequently for $\ds {\bK= \bigcup_{\{i_1,
\dotsc, i_s\} \in \{1, \dotsc, r\}} \{k_{i_1\dotsc i_s}\}}$ we have
\begin{equation}\label{indekscase5}
I_V(p,T)_H=
 \left\{\begin{array}{ll}
 (-1)^{j_0(V''(p),T)}  &    \text{ for } H=\sone,\\
  & \\
  (-1)^{j_0(V''(p),T)} \cdot j_k(V''(p),T)  &  \text{ for }  H=\bZ_k \text{ and } k \notin \bK.
\end{array} \right.
\end{equation}
\end{Remark}

\begin{Lemma}\label{Ik=degk}
Fix $p\in \{p_1,\ldots,p_q, \infty\}$ and assume that  $\ds \sigma(V''(p)) \cap \left\{ \frac{4k^2\pi^2}{T^2} : \
k \in \bN \right\} =$ $\ds = \left\{ \frac{4k_1^2\pi^2}{T^2}, \dotsc,\frac{4k_r^2\pi^2}{T^2}\right\}.$ If the
stationary solution $p \in \h^1_T$ is an isolated critical point of $\Phi_V$ then for $\displaystyle{k \notin \bK
= \bigcup_{\{i_1, \dotsc i_s\} \in \{1, \dotsc, r\}} \{k_{i_1\dotsc i_s}\}}$ we have
$$I_V(p,T)_{\bZ_k}= \dg_{\bZ_k}(\nabla \Phi_V, D_{\alpha_p}(\h^1_T,p)),$$ where
$\alpha$ is a conveniently chosen radius and $\Phi_V$ is a functional given by \eqref{NFUN}.
\end{Lemma}
\begin{proof}
Notice that if $p$ is an isolated critical point of $\Phi_V$, then functional $\Phi_V$ satisfies conditions
(F.1)-(F.6). Applying Lemma \ref{ogolny}, we conclude

$$\dg(\nabla \Phi_V,D_{\alpha_p}(\h^1_T,p))=$$ $$=\dg(\nabla \varphi_p,D_{\alpha_p}(V_p,p)) \star
\dg(-Id, D_{\alpha_p}(\h^1_{T,-},0)).$$ Since $ k \notin \bK,$  $\big(V_p\big)_{\bZ_k}=\emptyset.$ Hence by Remark
\ref{goodrem} we have $\dg_{\bZ_k}(\nabla \varphi_p,D_{\alpha_p}(V_p,p))=0.$ Therefore taking into account formula
\eqref{m} we obtain the following
$$
\dg_{\bZ_k}(\nabla \Phi_V,D_{\alpha_p}(\h^1_T,p))=$$ $$= \dg_{\sone}(\nabla \varphi_p,D_{\alpha_p}(V_p,p)) \cdot
\dg_{\bZ_k}(-Id, D_{\alpha_p}(\h^1_{T,-},0))$$ for $ k \notin \bK.$ But from Fact \ref{faktdegizo} we obtain
$$\dg(-Id,D_{\alpha_p}(\h^1_{T,-},0))=$$ $$=
((-1)^{j_0(V''(p),T)},(-1)^{j_0(V''(p),T)}\cdot j_1  (V''(p),T) ,\dotsc,(-1)^{j_0(V''(p),T)}\cdot
j_k(V''(p),T),\dotsc).$$ Taking into consideration that
$$\dg_{\sone}(\nabla \varphi_p,D_{\alpha_p}(V_p,p)) \cdot (-1)^{j_0(V''(p),T)}=\deg_B(-V',D_{\alpha_p}(\bR^n,p),0)$$
we complete the proof.
\end{proof}

\nt Combining the above considerations with Theorem \ref{istnienie} we can formulate its more general version.

\begin{Theorem}\label{istnienie2}
Suppose that
\begin{enumerate}
\item[(i)] $\big(V'\big)^{-1}(0)= \{p_1, \dotsc, p_q\}$,
\item[(ii)] $V'(x)= V''(\infty) \cdot x+o(\|x\|)$ as $\|x\| \rightarrow \infty.$
\end{enumerate}
If there exists $k \in \bN$ such that $I_V(p,T)_k$ is defined for all $p \in \{p_1, \dotsc, p_q,\infty\}$ and
 $\ds I_V(\infty,T)_k \neq$ $\ds \sum_{i=1}^qI_V(p_i,T)_k,$ then
there exists at least one non-stationary $T$-periodic solution of \eqref{NNSO}.  Moreover, if $p_0 \in \{p_1,
\dotsc, p_q,\infty\}$ is not an isolated  $T$-periodic solution of \eqref{NNSO}, then there are numbers
$k_1,\ldots,k_r \in \bN$ such that  $\ds \sigma(V''(p_0)) \cap \left\{\frac{4 k^2 \pi^2}{T^2} : k \in \bN
\right\}=\left\{\frac{4 k_1^2 \pi^2}{T^2},\ldots,\frac{4 k_r^2 \pi^2}{T^2} \right\} $ and the mi\-ni\-mal period
of any solution sufficiently close to $p_0$ equals $\ds \frac{T}{k_{i_1 \ldots i_s}},$ for some $\{i_1, \ldots,
i_s\} \subset \{1,\ldots,r\}.$
\end{Theorem}
\begin{proof} Assume  that stationary solutions $p_1, \dotsc, p_q$ and $\infty$ are isolated critical points of $\Phi_V$. If
they were not isolated, we would obtain an infinite sequence of critical points of $\Phi_V$, i.e. a sequence of
$T$-periodic non-stationary solutions of \eqref{NNSO}, which is our assertion. Suppose, contrary to our claim,
that $p_1, \dotsc, p_q$ are the only $T$-periodic solutions of \eqref{NNSO}. As in the proof of Theorem
\ref{istnienie} we obtain
$$\dg(\nabla \Phi_V, D_{\alpha_{\infty}}(\h^1_T,\infty)) =
\sum_{i=1}^q \dg (\nabla \Phi_V, D_{\alpha_{p_i}}(\h^1_T,p_i)).$$ But since $\Phi_V$ satisfies assumptions of the
Lemma \ref{Ik=degk}, we have
$$\dg_{\bZ_k}(\nabla \Phi_V, D_{\alpha_{\infty}}(\h^1_T,\infty))=I_V(\infty,T)_{\bZ_k} \neq$$
$$\neq \sum_{i=1}^q I_V(p_i,T)_{\bZ_k}= \sum_{i=1}^q \dg_{\bZ_k} (\nabla \Phi_V, D_{\alpha_{p_i}}(\h^1_T,p_i)).$$
Hence, $\ds \dg(\nabla \Phi,D_{\alpha_{\infty}}(\h^1_T,\infty)) \neq \sum_{i=1}^q\dg(\nabla
\Phi_V,D_{\alpha_{p_i}}(\h^1_T,p_i)),$ a contradiction.

\nt Fix $p_0 \in \{p_1,\dotsc, p_q,\infty\}$ and assume that $p_0$ is not isolated in $(\nabla \Phi_V)^{-1}(0).$
Then from Lemma \ref{redukcja} it follows that $\ds \sigma(V''(p_0)) \cap \left\{\frac{4 k^2 \pi^2}{T^2} : k \in
\bN \right\}=\left\{\frac{4 k_1^2 \pi^2}{T^2},\ldots,\frac{4 k_r^2 \pi^2}{T^2} \right\}.$ Hence from Corollary
\ref{wniosek} it follows that $\ds \ker \nabla^2 \Phi_V(p_0) \approx \bigoplus_{i=1}^r \bR
\big[\mu_{V''(p_0)}\big(\frac{4 k_i^2 \pi^2 }{T^2}\big),k\big].$ That is why the isotropy group $\sone_x$ of any
 element of $\ds x \in \ker \nabla^2 \Phi_V(p_0) \setminus \{0\}$ is equal to $\bZ_{k_{i_1 \ldots i_s}}$ for some $\{i_1, \ldots,
i_s\} \subset \{1,\ldots,r\}.$ It is clear that if $x \in \h^1_T$ and $\sone_x=\bZ_{k_{i_1 \ldots i_s}}$ then the
minimal period of $x$ is equal to $\ds \frac{T}{k_{i_1 \ldots i_s}}.$ The rest of the proof is a direct
consequence of Corollary \ref{cpi}.
\end{proof}

\begin{Remark}\label{wnistnienie2}
Notice that the above theorem one can formulate in the following e\-qui\-va\-lent way. Assume that for  $p \in
\{p_1,\ldots,p_q,\infty\}$  $$\sigma(V''(p)) \cap \left\{\frac{4 k^2 \pi^2}{T^2} : k \in \bN
\right\}=\left\{\frac{4 k_{1,p}^2 \pi^2}{T^2},\ldots,\frac{4 k_{r_p,p}^2 \pi^2}{T^2} \right\}$$ and define $\ds
\bK=\bigcup_{p \in \{p_1, \dotsc, p_q,\infty\}} \bigcup_{\{i_1, \dotsc i_s\} \in \{1, \dotsc, r_p\}}
\{\gcd\{k_{i_1,p},\ldots, k_{i_s,p}\}\}.$  If there exists $k \notin \bK$ such that
$$\ind(-V',\infty)  \cdot
j_k\left(V''(\infty),T\right) \neq \displaystyle{\sum_{i=1}^q \ind(-V',p_i) \cdot j_k\left(V''(p_i),T\right)}$$
then there exists at least one non-stationary $T$-periodic solution of \eqref{NNSO}.
\end{Remark}

\br \label{rwc} Notice that we cannot prove Theorem \ref{istnienie2} using the Leray-Schauder degree because it
vanishes. In fact, the assumption of Theorem \ref{istnienie2} can be rewritten in the following way
$$\nabla_{\sone}-\mathrm{deg}\bigg(\nabla \Phi ,D_{\alpha_{\infty}}\big(\h^1_T,\infty\big) \setminus
cl \bigg(\bigcup_{i=1}^p D_{\alpha_{p_i}}\big(\h^1_T,p_i\big) \bigg) \neq \Theta \in U(\sone).$$ From
  Theorem \ref{rawa} it follows that
$$\deg_{\rm LS}\bigg(\nabla \Phi,D_{\alpha_{\infty}}\big(\h^1_T,\infty\big) \setminus
cl \bigg(\bigcup_{i=1}^p D_{\alpha_{p_i}}\big(\h^1_T,p_i\big)\bigg),0\bigg)=$$
$$=\deg_{\rm LS}\bigg(\nabla \Phi^{\sone},\bigg(D_{\alpha_{\infty}}\big(\h^1_T,\infty\big) \setminus
cl \bigg(\bigcup_{i=1}^p D_{\alpha_{p_i}}\big(\h^1_T,p_i\big)\bigg)\bigg)^{\sone},0\bigg)=$$
$$=\deg_{\rm B} \bigg(-V',D_{\alpha_{\infty}}\big(\bR^n,\infty\big) \setminus
cl \bigg(\bigcup_{i=1}^p D_{\alpha_{p_i}}\big(\bR^n,p_i\big)\bigg),0\bigg)=0 \in \bZ $$ because
$\ds\big(V'\big)^{-1}(0) \cap \bigg(D_{\alpha_{\infty}}\big(\bR^n,\infty\big) \setminus cl \bigg(\bigcup_{i=1}^p
D_{\alpha_{p_i}}\big(\bR^n,p_i\big)\bigg)\bigg)=\emptyset.$
 \er

\numsubsec
\subsection{Continuation  of Periodic Solutions of Nonlinear Equation}
\label{persolcon} In this section we study continuation of non-stationary $T$-periodic solutions of the    family
of Newtonian systems of the form
\begin{equation}\label{rodzina}
(F)_{\lambda} \:\:\: \begin{cases}
    \ddot u =  -V'_{\lambda}(u)&  \\
    u(0)=u(T) & \\
    \dot u (0)=\dot u(T)
 \end{cases}
\end{equation}
where $V_{\lambda} \in C^2(\mathbb{R}^n, \mathbb{R}), \lambda \in \mathbb{R}.$ By $I_{V_{\lambda}}(p)$ we denote
the index defined in the previous section.

\begin{Theorem}\label{twkontynuacja}
Assume that
\begin{enumerate}
\item[(1)] $(V_0')^{-1} (0) = \{p_1, \dotsc , p_q\},$
\item[(2)] $V_0'(x) = V_0''(\infty) \cdot x + o(\|x\|)$ as $\|x\| \rightarrow  \infty$,
\item[(3)] there exists $k \in \bN$ such that $\ds I_{V_0}(\infty,T)_{\bZ_k}
\neq \sum_{i=1}^q I_{V_0}(p_i,T)_{\bZ_k}.$
\end{enumerate}
Then there exists an infinite sequence of non-stationary $T$-periodic solutions of $(F)_0$ converging to some $p
\in \{p_1, \dotsc , p_q,\infty\}$ or there exist closed, connected sets $\cC^{\pm}$ such that
\begin{equation*}
\begin{split}
&\cC^- \subset \big(\h^1_T \times (-\infty, 0]\big) \cap \big(\nabla
\Phi_{V_\lambda}\big)^{-1}(0),\\
&\cC^+ \subset \big(\h^1_T \times [0, +\infty)\big)\cap \big(\nabla \Phi_{V_\lambda}\big)^{-1}(0).
\end{split}
\end{equation*}
Moreover, for $\cC=\cC^{\pm}$
\begin{enumerate}
\item[(C1)]$\ds \cC \cap \bigg(\bigg(D_{\alpha_{\infty}}\big(\h^1_T,\infty\big) \setminus
cl \bigg(\bigcup_{i=1}^q D_{\alpha_{p_i}}\big(\h^1_T,p_i\big)\bigg)\bigg) \times \{0\}\bigg) \neq \emptyset,$
\item[(C2)] either $\cC$ is not bounded or $\cC \cap \{p_1, \dotsc,p_q\} \neq \emptyset.$
\end{enumerate}
Additionally, if $p_0 \in \{p_1, \dotsc, p_q,\infty\}$ is not an isolated  $T$-periodic solution of $(F)_0$, then
there are numbers $k_1,\ldots,k_r \in \bN$ such that  $$\ds \sigma(V''(p_0)) \cap \left\{\frac{4 k^2 \pi^2}{T^2} :
k \in \bN \right\}= \left\{\frac{4 k_1^2 \pi^2}{T^2},\ldots,\frac{4 k_r^2 \pi^2}{T^2} \right\}$$ and the
mi\-ni\-mal period of any solution sufficiently close to $p_0$ equals $\ds \frac{T}{k_{i_1 \ldots i_s}},$ for some
$\{i_1, \ldots, i_s\} \subset \{1,\ldots,r\}.$
\end{Theorem}
\begin{proof}
Consider functional $\Phi_{V_{\lambda}} \in C^2_{\sone}(\h^1_T , \mathbb{R})$ given by \eqref{NFUN} and suppose
that stationary solutions $p_1,\dotsc, p_q$ and $\infty$ are isolated solutions of \eqref{rodzina} on level
$\lambda_0=0,$ i.e. they are isolated critical points of $\Phi_{V_0}.$ Thus we can choose
$\alpha_{\infty},\alpha_{p_i} > 0, i=1,\dotsc,q$ such that
\begin{enumerate}
\item[(i)] $\big(\nabla \Phi_{V_0}\big)^{-1}(0) \cap(\h^1_T \setminus
D_{\alpha_{\infty}}(\h^1_T,\infty))= \emptyset,$
\item[(ii)] $\big(\nabla \Phi_{V_0}\big)^{-1}(0) \cap D_{\alpha_{p_i}}(\h^1_T,p_i)
= \{p_i\}, i= 1,\dotsc,q$,
\item[(iii)] $D_{\alpha_{p_i}}(\h^1_T,p_i) \cap D_{\alpha_{p_j}}(\h^1_T,p_j)=
\emptyset$, for $i \neq j.$
\end{enumerate}
Set $\Omega = D_{\alpha_{\infty}}(\h^1_T,\infty) \setminus cl\bigg(\displaystyle{\bigcup_{i=1}^q
D_{\alpha_{p_i}}(\h^1_T,p_i)}\bigg)$. Notice that from the assumption (3) and Facts \ref{I=deg}, \ref{Ik=degk} we
obtain $\ds\dg(\nabla\Phi_{V_0},\Omega) \neq \Theta \in U(\sone).$ The rest of the proof is a direct consequence
of Theorem \ref{abscont}. The second part of the proof is in fact the same as the proof of Theorem
\ref{istnienie2}.
\end{proof}

\begin{Corollary} \label{ghtr} Let assumptions of Theorem \ref{twkontynuacja}
be satisfied. If moreover $\ds \sigma(V''(p)) \cap \left\{ \frac{4k^2\pi^2}{T^2} : \ k \in \bN \right\}=
\emptyset$ for any $p\in \{p_1,\dotsc,p_q,\infty\}$ then there exist closed connected sets $\cC^- \subset (\h^1_T
\times (-\infty, 0]) \cap \big(\nabla \Phi_{V_\lambda}\big)^{-1}(0), \cC^+ \subset (\h^1_T \times [0,
+\infty))\cap \big(\nabla \Phi_{V_\lambda}\big)^{-1}(0)$ with properties (C1), (C2).
\end{Corollary}

\begin{proof} From Lemma \ref{redukcja} we conclude that any $p\in \{p_1,\dotsc,p_q, \infty\}$
is an isolated critical point of $\Phi_{V_0}$ i.e. an isolated solution of $\ddot{u}=-V'_{0}(u)$.
\end{proof}

\begin{Remark}
If $\cC=\cC^{\pm}$ in Theorem \ref{twkontynuacja} is bounded, then symmetry breaking phenomenon occurs, i.e. $\cC$
contains solutions with different minimal periods. Indeed, in this case from (C2) we obtain that $\cC$ contains
stationary solution (whose isotropy group in $\h^1_T$ is equal to $\sone$) and from (C1) - non-stationary solution
(with isotropy group $\bZ_k$ for some $k \in \bN$, which means that its minimal period is equal to $\frac{T}{k}$).
\end{Remark}

\br Similarly as in the previous section the Leray-Schauder degree is not applicable in the proof of Theorem
\ref{twkontynuacja}.\er

\numsec

\section{Illustration}
\label{illu}

In this section we   illustrate abstract results proved in the previous parts of this article. Namely, we study
non-stationary $T$-periodic solutions of the following system
\begin{equation}\label{sitnikov}
 \begin{cases}
    \ddot u =  -V' (u)&  \\
    u(0)=u(T) & \\
    \dot u (0)=\dot u(T)
 \end{cases}
\end{equation}
where  potential  $ V: \bR^n \rightarrow \bR$ is defined as follows
\begin{equation}\label{sitnikovpotential}
V (x)= \frac{1}{2} <V''(\infty) x,x> + W (x) =\frac{1}{2} <V''(\infty) x,x> + \frac{-1}{\sqrt{\|x\|^2 +a}}
\end{equation}
where $a > 0$ and  $V''(\infty)$ is a real symmetric $(n \times n)$-matrix.

\nt In the following four lemmas we study properties of the functional $V.$

\bl \label{exprop1} If potential $V$ is given by \eqref{sitnikovpotential} then $$\ds V'(x)= V''(\infty) \cdot x +
W'(x)  = V''(\infty) \cdot x + \frac{1}{(\|x\|^2 +a)^{3/2}} \cdot x = V''(\infty) \cdot x + o(\|x\|) \text{ as }
\|x\| \rightarrow \infty.$$ Additionally, there is an orthogonal matrix $P \in O(n,\bR)$ such that
$$ V'(x)= P \cdot \bigg( J(V''(\infty)) + \frac{1}{(\|x\|^2 +a)^{3/2}}\ \cdot Id_{\bR^n}\bigg) \cdot P^{-1} \cdot
  x,$$ where   $J(V''(\infty))=\diag\{\lambda_1,\ldots,\lambda_n\}.$
Moreover, we have
\begin{enumerate}
  \item $\ds V''(x)= V''(\infty) + W''(x),$
  \item $\ds \big(W\big)''_{x_i x_i}(x)={\frac {-3{x_{{i}}}^{2}}{ \left(\|x\|^2+a\right)^{5/2}}}+
 \frac{1}{\left(\|x\|^2+a \right)^{3/2}}$ for $i=1,\ldots,n,$
 \item $\ds \big(W\big)''_{x_i x_j}(x)={\frac {-3 x_{{i}}x_{{j}}}{ \left(\|x\|^2+a\right)^{5/2}}} $
 for $i \neq j, i,j = 1,\ldots,n.$
\end{enumerate}
\el

\nt The easy proof of the above lemma is left to the reader.

\bl \label{exprop2}  Assume that $\mu_{V''(\infty)}(\lambda_i)=1$ for any $\ds \lambda_i \in \sigma(V''(\infty))
\cap \left[- \frac{1}{\sqrt{a^3}}, 0\right).$ Then $\# (V')^{-1}(0)  < \infty.$ Moreover, if $x \in (V')^{-1}(0)$
then $x=Py$ where for $i=1,\ldots,n$ we have
$$y_i =
  \begin{cases}
    0, & \text{ if } \lambda_i \geq 0 \text{ or }  \lambda_i <\ds  - \frac{1}{\sqrt{a^3}}, \\
     &   \\
    0 \text{ or } \pm \sqrt{\frac{\ds1}{\ds \sqrt[3]{\lambda_i^2}}-a}, & \text{otherwise}.
  \end{cases}
$$
\el
\begin{proof} Let $P \in O(n,\bR)$ be as in Lemma \ref{exprop1}.

\nt Hence
  $V'(y)=0$ iff $P^{-1} V'(Py)=0$ and moreover
$$P^{-1} V'(Py)=
  \bigg( J(V''(\infty)) + \frac{1}{(\|y\|^2 +a)^{3/2}}\ \cdot Id_{\bR^n}\bigg) \cdot y.$$

\nt Fix $\lambda_i \in \sigma_+(V''(\infty)) \cup \{0\}$ and notice that $\ds   \lambda_i +  \frac{1}{(\|y\|^2
+a)^{3/2}} > 0$ for any $y \in \bR^n.$ Assume now that $\lambda_i \in \sigma_-(V''(\infty)).$ It is clear that $$
\lambda_i + \frac{1}{(\|y\|^2 +a)^{3/2}} = 0 \:\: \mathrm{ iff } \: \|y\|^2 = \frac{1}{\ds
\sqrt[3]{\lambda_i^2}}-a$$ and consequently if $\ds a > \frac{1}{\ds\sqrt[3]{\lambda_i^2}}$ then $\ds \lambda_i +
\frac{1}{(\|y\|^2 +a)^{3/2}} < 0.$

\nt Taking into account that $\mu_{V''(\infty)}(\lambda_i)=1$ for any $\lambda_i \in \sigma_-(V''(\infty))$ and
the above we obtain that
$$\bigg( J(V''(\infty)) + \frac{1}{(\|y\|^2 +a)^{3/2}} \cdot Id_{\bR^n}\bigg)\left[\begin{array}{c}
  y_1 \\
  \vdots \\
  y_n
\end{array} \right]=\left[\begin{array}{c}
  0 \\
  \vdots \\
  0
\end{array} \right]$$
if and only if for $i=1,\ldots,n$ we have $$y_i =
  \begin{cases}
    0, & \text{ if } \lambda_i \geq 0 \text{ or }  \lambda_i <\ds  - \frac{1}{\sqrt{a^3}}, \\
     &   \\
    0 \text{ or } \pm \sqrt{\frac{1}{\ds \sqrt[3]{\lambda_i^2}}-a}, & \text{otherwise},
  \end{cases}
$$
which completes the proof.
\end{proof}

\bl \label{exprop3}  If $\ds \sigma(V''(\infty)) \cap \left[- \frac{1}{\sqrt{a^3}}, 0\right)=\emptyset$ then
$(V')^{-1}(0)=\{\Theta\}.$\el
\begin{proof}   By Lemma \ref{exprop4} we have
\begin{equation}\label{wert}
  V'(x)= P \cdot \bigg( J(V''(\infty)) + \frac{1}{(\|x\|^2 +a)^{3/2}}\ \cdot Id_{\bR^n}\bigg) \cdot P^{-1}
\cdot x
\end{equation}
Since $\ds \sigma(V''(\infty)) \cap \left[- \frac{1}{\sqrt{a^3}}, 0\right)=\emptyset,$ $\ds \lambda_i +
\frac{1}{(\|x\|^2 +a)^{3/2}}  \neq 0$ for $i=1,\ldots,n$ and any $x \in \bR^n.$ Hence matrix $\ds P \cdot \bigg(
J(V''(\infty)) + \frac{1}{(\|x\|^2 +a)^{3/2}}\ \cdot Id_{\bR^n}\bigg) \cdot P^{-1}$ is nondegenerate for any $x
\in \bR^n.$ Taking into account \eqref{wert} we complete the proof.
\end{proof}

\bl \label{exprop4} Under the above assumptions $\ds \ind(-V',\infty)=(-1)^{n-m^-(V''(\infty))}.$\el
\begin{proof} Notice that for   sufficiently large $\|x\|$ and $\lambda_i \in \sigma(V''(\infty)) \setminus \{0\}$ we have
$$\sign \bigg(\lambda_i +  \frac{\ds1}{\ds(\|x\|^2 +a)^{3/2}}\bigg) = \sign \lambda_i.$$

\nt For $i=1,\ldots,n$ define $\psi_i :[0,1] \rightarrow \bR \setminus \{0\}$ as follows
$$\psi_i(t)=
  \begin{cases}
  t \cdot   \sign \lambda_i + (1-t) \cdot \bigg(\lambda_i +  \frac{\ds1}{\ds(\|x\|^2 +a)^{3/2}}\bigg) & \text{ if } \lambda_i \neq 0, \\
  t  +     \frac{\ds 1 -t}{\ds(\|x\|^2 +a)^{3/2}}     & \text{ if } \lambda_i =0.
  \end{cases}
$$
Since the above, for any $t \in [0,1]$ and sufficiently large $\|x\|$ we have $$\det \left( P \cdot
\left[\begin{array}{ccc}
  \psi_1(t)& \ldots & 0 \\
  \vdots & \ddots & \vdots \\
  0 & \ldots & \psi_n(t)
\end{array}\right] \cdot P^{-1} \right) = \prod_{i=1}^n \psi_i(t) \neq 0.$$
Define  a map $\Psi : \bR^n \times [0,1] \rightarrow \bR^n$ as follows
$$\Psi(x,t)= \left(P \cdot
\left[\begin{array}{ccc}
  \psi_1(t)& \ldots & 0 \\
  \vdots & \ddots & \vdots \\
  0 & \ldots & \psi_n(t)
\end{array}\right] \cdot P^{-1} \right) \cdot x.$$
It is easy to verify that
\begin{enumerate}
  \item $\Psi(\cdot,0)=V'(\cdot),$
  \item there is $\beta_0 > 0$ such that for any $t \in [0,1]$ and any $x \in \bR^n$ such that $\|x\| > \beta_0$
  $\Psi(x,t) \neq 0,$
  \item $\ind (-V',\infty)=\ind(-\Psi(\cdot,0),\infty)=\ind(-\Psi(\cdot,1),\infty)=(-1)^{n-m^-(V''(\infty))},$
\end{enumerate}
which completes the proof.
\end{proof}

\bex \label{ist1} In this example we study system \eqref{sitnikov} with resonance at the infinity ($V''(\infty)$
is degenerate). Potential $V$ is a Morse function i.e. all the critical points of $V$ are nondegenerate. The
origin $\Theta \in \bR^n,$ treated as a constant function, is a resonant stationary solution of \eqref{sitnikov}
i. e. $\sigma (V''(\Theta)) \cap \{k^2 : k \in \bN\} \neq \emptyset.$ Consider system \eqref{sitnikov} with $n=4,
a=1, T=2\pi$ and
$$V''(\infty)=\left[\begin{array}{cccc}\ds
 \frac{7}{2} & 0 & 0 & 0 \\
  0 & -2 & 0 & 0 \\
  0 & 0 & 0 & 0 \\
  0 & 0 & 0 & -\frac{\ds 1}{\ds 2\sqrt{2}}
\end{array} \right].$$  From Lemma \ref{exprop2} it follows that
$(V')^{-1}(0)=\{\Theta, \pm e_4\}.$ Moreover, by  Lemma \ref{exprop1} we have
$$V''(\Theta)=V''(\infty)+W''(\Theta)=V''(\infty) + Id=\left[\begin{array}{cccc}\ds
 \frac{9}{2} & 0 & 0 & 0 \\
  0 & -1 & 0 & 0 \\
  0 & 0 & 1 & 0 \\
  0 & 0 & 0 & 1-\frac{\ds 1}{\ds 2\sqrt{2}}
\end{array} \right]$$ and
$$V''(\pm e_4)=V''(\infty)+W''(\pm e_4)= \left[\begin{array}{cccc}\ds
 \frac{7}{2} + \frac{\ds 1}{\ds 2\sqrt{2}}& 0 & 0 & 0 \\
  0 & -2+\frac{\ds 1}{\ds 2\sqrt{2}} & 0 & 0 \\
  0 & 0 & \frac{\ds 1}{\ds 2\sqrt{2}} & 0 \\
  0 & 0 & 0 & \frac{\ds -3}{\ds 4 \sqrt{2}}
\end{array} \right].$$
Since $\det(-V''(\Theta)) < 0, \ind(-V',\Theta)=-1.$ Additionally condition $\det (-V''(\pm e_4)) > 0$ implies
$\ind(-V',\pm e_4)=1.$ Moreover, from Lemma \ref{exprop4} it follows that $\ind(-V',\infty)=1.$ Finally notice
that
$$j_k(V''(\Theta),2\pi)=
  \begin{cases}
     1 & \text{ if } k=1,2 \\
     0 & \text{otherwise,}
  \end{cases}, \:\:
j_k(V''(\infty),2\pi)= j_k(V''(\pm e_4),2\pi)=
  \begin{cases}
     1 & \text{ if } k=1, \\
     0 & \text{otherwise.}
  \end{cases}
$$
It is clear that $\sigma(V''(\Theta)) \cap \{k^2 : k \in \bN\}=\{1\}$ and that $$\sigma(V''(\pm e_4)) \cap \{k^2 :
k \in \bN\}= \sigma(V''(\infty)) \cap \{k^2 : k \in \bN\}=\emptyset$$ Moreover,
$$\ind(-V',\infty) \cdot j_2(V''(\infty),2\pi)= 1 \cdot 0 = 0\neq  -1= (-1) \cdot 1 + 1 \cdot 0 + 1 \cdot 0 =$$
$$=\ind(-V',\Theta) \cdot j_2(V''(\Theta),2\pi) +
\ind(-V',e_4) \cdot j_2(V''(e_4),2\pi) + $$ $$+\ind(-V',-e_4) \cdot j_2(V''(-e_4),2\pi).$$ We have just shown that
all the assumption of Theorem  \ref{istnienie2}  are satisfied. Therefore there is at least one non-stationary $2
\pi $-periodic solution of system \eqref{sitnikov}.  \eex

\bex  \label{ist2} In this example we study system \eqref{sitnikov} with resonance at the infinity ($V''(\infty)$
is degenerate). Potential $V$ is not a Morse function because $\Theta \in \bR^4$ is a degenerate critical point of
$V.$ The origin $\Theta \in \bR^n,$ treated as a constant function, is a resonant stationary solution of
\eqref{sitnikov} i. e. $\sigma (V''(\Theta)) \cap \{k^2 : k \in \bN\} \neq \emptyset.$ Consider system
\eqref{sitnikov} with $n=4, a=1, T=2\pi$ and $$ V''(\infty)=\left[\begin{array}{cccc}\ds
 \frac{7}{2} & 0 & 0 & 0 \\
  0 & -1 & 0 & 0 \\
  0 & 0 & 0 & 0 \\
  0 & 0 & 0 & -\frac{\ds 1}{\ds 2\sqrt{2}}
\end{array} \right].$$  From Lemma \ref{exprop1} it follows that
$(V')^{-1}(0)=\{\Theta, \pm e_4\}.$ Moreover, by  Lemma \ref{exprop1} we have
$$V''(\Theta)=V''(\infty)+W''(\Theta)=V''(\infty) + Id=\left[\begin{array}{cccc}\ds
 \frac{9}{2} & 0 & 0 & 0 \\
  0 & 0 & 0 & 0 \\
  0 & 0 & 1 & 0 \\
  0 & 0 & 0 & 1-\frac{\ds 1}{\ds 2\sqrt{2}}
\end{array} \right]$$ and
$$V''(\pm e_4)=V''(\infty)+W''(\pm e_4)= \left[\begin{array}{cccc}\ds
 \frac{7}{2} + \frac{\ds 1}{\ds 2\sqrt{2}}& 0 & 0 & 0 \\
  0 & -1+\frac{\ds 1}{\ds 2\sqrt{2}} & 0 & 0 \\
  0 & 0 & \frac{\ds 1}{\ds 2\sqrt{2}} & 0 \\
  0 & 0 & 0 & \frac{\ds -3}{\ds 4 \sqrt{2}}
\end{array} \right].$$ Since $\det ( -V''(\pm e_4)) > 0,$ $\ind(-V',\pm e_4)=1.$
By Lemma \ref{exprop4} we have $\ind(-V',\infty)=1.$
Consequently by Remark \ref{sumfo} we obtain that $\ind(-V',\Theta)=-1.$ Finally notice that
$$j_k(V''(\Theta),2\pi)=
  \begin{cases}
     1 & \text{ if } k=1,2 \\
     0 & \text{otherwise,}
  \end{cases}, \:\:
j_k(V''(\infty),2\pi)= j_k(V''(\pm e_4),2\pi)=
  \begin{cases}
     1 & \text{ if } k=1, \\
     0 & \text{otherwise.}
  \end{cases}
$$
It is evident that $\ds \sigma(V''(\Theta)) \cap \{k^2 : k \in \bN\}=\{1\}, \det (V''(\Theta))=0 $ and that
$\sigma(V''(\pm e_4)) \cap $ $\{k^2 : k \in \bN\}=$ $ \sigma(V''(\infty)) \cap \{k^2 : k \in \bN\}=\emptyset.$
Moreover,
$$\ind(-V',\infty) \cdot j_2(V''(\infty),2\pi) = 1 \cdot 0 = 0\neq -1= (-1) \cdot 1 + 1 \cdot 0 + 1 \cdot 0 =$$
$$=\ind(-V',\Theta) \cdot j_2(V''(\Theta),2\pi) +
\ind(-V',e_4) \cdot j_2(V''(e_4),2\pi) +$$ $$+ \ind(-V',-e_4) \cdot j_2(V''(-e_4),2\pi).$$ We have just shown that
all the assumption of Theorem  \ref{istnienie2}  are satisfied. Therefore there is at least one non-stationary $2
\pi $-periodic solution of system \eqref{sitnikov}. \eex

\bex  \label{ist3} Let us put in \eqref{sitnikov}  $\ds n=1, V''(\infty)=0, a=\frac{1}{4}.$ Under these
assumptions \eqref{sitnikov}   is the Sitnikov circular problem studied for example in \cite{[COLL]},
\cite{[LLIB-ORT]}. It is easy to check that $(V')^{-1}(0)=\{0\}, V''(0)=8, \ind(-V',0)=\ind(-V',\infty)=-1.$
Finally notice that
$$j_k(V''(0),T)= j_k(8,T)=
  \begin{cases}
     1 & \text{ if } k < \frac{\ds T \sqrt{2}}{\ds \pi} \\
     0 & \text{otherwise,}
  \end{cases}, $$ and $
j_k(V''(\infty),T)=j_k(0,T)=0 $ for any $k \in \bN.$ Summing up, if  $\ds T > \frac{\pi}{\sqrt{2}}$ then
$$\ind(-V',0) \cdot j_1(V''(0),T)= (-1) \cdot 1 =-1 \neq 0=(-1) \cdot 0 = \ind(-V',\infty) \cdot j_1(V''(\infty),T).$$

\nt We have just shown that all the assumption of Theorem  \ref{istnienie2}  are satisfied. Therefore there is at
least one non-stationary $T$-periodic solution of the circular Sitnikov problem for any $\ds T >
\frac{\pi}{\sqrt{2}}$. \eex

\nt In th rest of this  section we study continuation of non-stationary $T$-periodic solutions of the  family of
Newtonian systems of the form
\begin{equation}\label{sitnikovfamily}
 \begin{cases}
    \ddot u =  -V'_{\lambda}(u)&  \\
    u(0)=u(T) & \\
    \dot u (0)=\dot u(T)
 \end{cases}
\end{equation}
where $V_{\lambda} \in C^2(\mathbb{R}^n, \mathbb{R}), \lambda \in \mathbb{R}$ and $V_0$ is given by formula
\eqref{sitnikovpotential}.

\bex  \label{cont1} Let us consider system \eqref{sitnikovfamily} with $V_0$ satisfying all the assumptions of
po\-ten\-tial $V$ considered in Example \ref{ist1} or in Example \ref{ist2}. It is easy to verify that system
\eqref{sitnikovfamily} satisfies all the assumptions of Theorem  \ref{twkontynuacja}. That is why the set of
non-stationary $2\pi$-periodic solutions of system \eqref{sitnikovfamily} satisfies the alternative from  of
Theorem \ref{twkontynuacja}.  \eex

\bex \label{cont2} Let us consider system \eqref{sitnikovfamily} with $V_0$ satisfying all the assumptions of the
Sitnikov  po\-ten\-tial $V$ considered in Example \ref{ist3}. It is easy to verify that system
\eqref{sitnikovfamily} satisfies all the assumptions of Theorem  \ref{twkontynuacja} and Corollary \ref{ghtr}.
Hence for any $\ds T > \frac{\pi}{\sqrt{2}}$  there are  closed connected sets $\cC^{\pm}$ of non-stationary
$T$-periodic solutions of system \eqref{sitnikovfamily} with properties $(C1), (C2)$ from Theorem
\ref{twkontynuacja}.\eex

\end{document}